\def\todaysdate{28\textsuperscript{th} November 2026}
\documentclass[reqno,a4paper]{article}
\usepackage{etex}
\usepackage[T1]{fontenc}
\usepackage[utf8]{inputenc}
\usepackage{lmodern}

\usepackage[style=alphabetic,giveninits=true,backref,backend=biber,isbn=false,url=false,maxbibnames=99,maxalphanames=4]{biblatex}

\renewbibmacro{in:}{}
\newbibmacro{string+doi}[1]{\iffieldundef{doi}{\iffieldundef{url}{#1}{\href{\thefield{url}}{#1}}}{\href{http://dx.doi.org/\thefield{doi}}{#1}}}
\DeclareFieldFormat*{title}{\usebibmacro{string+doi}{\emph{#1}}}
\DefineBibliographyStrings{english}{%
  backrefpage = {$\uparrow$},
  backrefpages = {$\uparrow$},
}

\usepackage{amssymb,amsmath,amsthm}
\usepackage{thmtools,xcolor}
\usepackage{units}
\usepackage{graphicx}
\usepackage[all]{xy}
\usepackage{tikz}
\usetikzlibrary{arrows,calc,positioning,decorations.pathreplacing}
\usepackage{tikz-cd}
\usepackage{relsize}
\usepackage{mhsetup}
\usepackage{mathtools}
\usepackage{stmaryrd}
\usepackage{tocloft}
\usepackage{paralist} % For compactitem environment
\usepackage[left=3cm,right=3cm,top=3cm,bottom=3cm]{geometry} % Decrease the margins %
\usepackage[bottom]{footmisc}
\usepackage[colorlinks,citecolor=blue,linkcolor=blue,urlcolor=blue,filecolor=blue,breaklinks]{hyperref}
\usepackage{footnotebackref}
\usepackage[format=plain,indention=1em,labelsep=quad,labelfont={up},textfont={footnotesize},margin=2em]{caption}
\usepackage{subcaption}
\usepackage[mathscr]{euscript}
\usepackage{accents}
\usepackage{xparse}
\usepackage{nicematrix}
\usepackage{marginnote}

\definecolor{lightblue}{rgb}{0.8,0.8,1}

\numberwithin{equation}{section}
\numberwithin{figure}{section}

\definecolor{vdarkred}{rgb}{0.7,0,0}
\declaretheoremstyle[
  spaceabove=\topsep,
  spacebelow=\topsep,
  headpunct=,
  numbered=no,
  postheadspace=1ex,
  headfont=\color{vdarkred}\normalfont\bfseries,
  bodyfont=\normalfont\itshape,
]{colored}
\declaretheoremstyle[
  spaceabove=\topsep,
  spacebelow=\topsep,
  headpunct=,
  numbered=no,
  postheadspace=1ex,
  headfont=\normalfont\bfseries,
  bodyfont=\normalfont\itshape,
]{italic}
\declaretheoremstyle[
  spaceabove=\topsep,
  spacebelow=\topsep,
  headpunct=,
  numbered=no,
  postheadspace=1ex,
  headfont=\normalfont\bfseries,
  bodyfont=\normalfont\upshape,
]{upright}

\declaretheorem[style=italic,name=Theorem,numbered=yes,numberwithin=section]{thm}
\declaretheorem[style=italic,name=Lemma,numbered=yes,numberlike=thm]{lem}
\declaretheorem[style=italic,name=Proposition,numbered=yes,numberlike=thm]{prop}
\declaretheorem[style=italic,name=Corollary,numbered=yes,numberlike=thm]{coro}

\declaretheorem[style=italic,name=Theorem,numbered=yes,numberwithin=section]{athm}

\declaretheorem[style=upright,name=Remark,numbered=yes,numberlike=thm]{rmk}

\declaretheorem[style=upright,name=Notation,numbered=yes,numberlike=thm]{notation}

\makeatletter
\renewcommand*{\@seccntformat}[1]{\upshape\csname the#1\endcsname.\hspace{1ex}}
\renewcommand*{\section}{\@startsection{section}{1}{\z@}%
	{2.5ex \@plus 1ex \@minus 0.2ex}%
	{1.5ex \@plus 0.2ex}%
	{\normalfont\Large\bfseries}}
\renewcommand*{\subsection}{\@startsection{subsection}{2}{\z@}%
	{2.5ex \@plus 1ex \@minus 0.2ex}%
	{1.5ex \@plus 0.2ex}%
	{\normalfont\large\bfseries}}
\renewcommand*{\subsubsection}{\@startsection{subsubsection}{3}{\z@}%
	{2.5ex \@plus 1ex \@minus 0.2ex}%
	{1.5ex \@plus 0.2ex}%
	{\normalfont\normalsize\bfseries}}
\newcommand*{\subsubsubsection}{\@startsection{paragraph}{4}{\z@}%
	{2.5ex \@plus 1ex \@minus 0.2ex}%
	{1.5ex \@plus 0.2ex}%
	{\normalfont\normalsize\bfseries}}
\newcommand*{\subsubsubsubsection}{\@startsection{paragraph}{5}{\z@}%
	{2.5ex \@plus 1ex \@minus 0.2ex}%
	{1.5ex \@plus 0.2ex}%
	{\normalfont\normalsize\bfseries}}
\makeatother

\newcommand{\Aut}{\mathrm{Aut}}

\newcommand{\B}{\mathbf{B}}

\newcommand{\Sym}{\mathfrak{S}}
\newcommand{\Hom}{\mathrm{Hom}}

\newcommand{\St}{\mathrm{St}}
\newcommand{\aDiag}{\overline{\Delta}}

\newcommand{\Image}{\mathrm{im}}

\newcommand{\incl}[3][right]%
{%
\draw[<-,>=#1 hook] #2 to ($ #2!0.5!#3 $);
\draw[->,>=stealth'] ($ #2!0.5!#3 $) to #3;%
}
\newcommand{\inclusion}[5][right]%
{%
\draw[<-,>=#1 hook] #4 to ($ #4!0.5!#5 $) node[#2,font=\small]{#3};
\draw[->,>=stealth'] ($ #4!0.5!#5 $) to #5;%
}

\newcommand{\longtwoheadrightarrow}{%
  \longrightarrow\mathrel{\mkern-14mu}\to
}
\tikzset{
  Isom/.style={
    draw=none,
    every to/.append style={
      edge node={node [sloped, allow upside down, auto=false]{$\cong$}}}
  }
}

\newenvironment{itemizeb}%
{\begin{compactitem}

}%
{\end{compactitem}}
{\begin{compactitem}[#1]

}%
{\end{compactitem}}
{\begin{compactdesc}

}%
{\end{compactdesc}}

\newcommand{\cP}{\mathcal{P}}

\newcommand{\bC}{\mathbb{C}}

\newcommand{\bQ}{\mathbb{Q}}

\newcommand{\bZ}{\mathbb{Z}}

\newcommand{\LKB}{\mathfrak{LKB}}
\newcommand{\Bur}{\mathfrak{Bur}}
\newcommand{\UBur}{\mathfrak{UBur}}
\newcommand{\TYM}{\mathfrak{TYM}}
\newcommand{\tr}{\mathrm{tr}}
\newcommand{\inv}{\mathrm{inv}}

\newcommand{\bzeta}{\boldsymbol{\zeta}}

\renewcommand{\geq}{\geqslant}
\renewcommand{\leq}{\leqslant}
\renewcommand{\footnoterule}{%
  \kern -3pt
  \hrule width \textwidth height 0.4pt
  \kern 2.6pt
}

\NewDocumentCommand{\covlr}{O{\bullet} O{\bullet}}{%
\ensuremath{{}_{#1}\mathrm{Cov}_{#2}}%
}
\NewDocumentCommand{\covl}{O{\bullet}}{%
\ensuremath{{}_{#1}\mathrm{Cov}}%
}
\NewDocumentCommand{\covr}{O{\bullet}}{%
\ensuremath{\mathrm{Cov}_{#1}}%
}
\NewDocumentCommand{\toplr}{O{\bullet} O{\bullet}}{%
\ensuremath{{}_{#1}\mathrm{Top}_{#2}}%
}
\NewDocumentCommand{\topl}{O{\bullet}}{%
\ensuremath{{}_{#1}\mathrm{Top}}%
}
\NewDocumentCommand{\topr}{O{\bullet}}{%
\ensuremath{\mathrm{Top}_{#1}}%
}
\NewDocumentCommand{\modlr}{O{\bullet} O{\bullet}}{%
\ensuremath{{}_{#1}\mathrm{Mod}_{#2}}%
}
\NewDocumentCommand{\modl}{O{\bullet}}{%
\ensuremath{{}_{#1}\mathrm{Mod}}%
}
\NewDocumentCommand{\modr}{O{\bullet}}{%
\ensuremath{\mathrm{Mod}_{#1}}%
}

\NewDocumentCommand{\lmod}{O{\bullet}}{%
\ensuremath{{}\textrm{-}\mathrm{Mod}}%
}

\NewDocumentCommand{\topmodlr}{O{\bullet} O{\bullet} O{\bullet}}{%
\ensuremath{{}_{#1}\mathrm{Top}_{#2}\mathrm{Mod}_{#3}}%
}
\NewDocumentCommand{\topmodl}{O{\bullet} O{\bullet}}{%
\ensuremath{{}_{#1}\mathrm{Top}_{#2}\mathrm{Mod}}%
}
\NewDocumentCommand{\topmodr}{O{\bullet} O{\bullet}}{%
\ensuremath{\mathrm{Top}_{#1}\mathrm{Mod}_{#2}}%
}
\NewDocumentCommand{\setlr}{O{\bullet} O{\bullet}}{%
\ensuremath{{}_{#1}\mathrm{Set}_{#2}}%
}
\NewDocumentCommand{\setl}{O{\bullet}}{%
\ensuremath{{}_{#1}\mathrm{Set}}%
}
\NewDocumentCommand{\setr}{O{\bullet}}{%
\ensuremath{\mathrm{Set}_{#1}}%
}
\definecolor{dblue}{rgb}{0.4,0,1}

\definecolor{dgreen}{rgb}{0,0.5,0}

\begin{document}
\title{\LARGE\bfseries Invariants and coinvariants\\ of braid representations}
\author{\normalsize Anh Trong Nam Hoang and Arthur Souli{\'e}}
\date{\normalsize\todaysdate}
\maketitle
\makeatletter
\renewcommand*{\BHFN@OldMakefntext}{}
\makeatother

\begin{abstract}
We compute the invariants and coinvariants of several families of braid group representations over their natural coefficient rings, as well as those of their specialisations at complex values of the parameters. These include the Tong--Yang--Ma, Burau, and Lawrence--Krammer--Bigelow representations, together with their duals and various tensor products. We identify the parameter values at which the dimensions of the invariant and coinvariant spaces change. As applications, our computations provide reducibility results for certain specialisations and exhibit families of coefficients for which twisted homological stability fails already in degree zero.
\end{abstract}

\tableofcontents

%-------------------------------------------------------

\section{Introduction}\label{sec:intro}

The {\it braid group} $\B_{n}$ admits the presentation:
\[\B_{n} = \langle \sigma_{1}, \dots, \sigma_{n-1} \mid \sigma_{i} \sigma_{j} = \sigma_{j} \sigma_{i} \text{ if } |i-j| > 1; \sigma_{i} \sigma_{i+1} \sigma_{i} = \sigma_{i+1} \sigma_{i} \sigma_{i+1} \rangle.\]
Recall that a representation of $\B_{n}$ over a unital ring $R$ is a group homomorphism $\rho \colon \B_{n} \to \Aut_{R}(V)$ where $V$ is an $R$-module. Such a representation is said to be \emph{irreducible} if $V$ has no proper nonzero $\B_{n}$-submodule. Understanding irreducible representations, and more generally the subrepresentations and quotients of a given representation, is a fundamental problem in the representation theory of braid groups.
 
The representation theory of braid groups is a rich and extensively studied subject; see, for instance, the survey by Birman and Brendle~\cite[\S 4]{BirmanBrendlesurvey} and the monograph by Kassel and Turaev~\cite{KasselTuraev}. Linear representations provide powerful tools for investigating braid groups, allowing group-theoretic questions to be approached through linear algebra. Their classification is, however, a difficult problem: for $n\geq 6$, the representation theory of $\B_{n}$ is known to be \emph{wild}, by a result of Erdmann and Nakano~\cite{ErdmannNakano}, in the sense of~\cite{BelitskiiSergeichuk}. Nevertheless, the classification of \emph{low-dimensional} representations of $\B_{n}$, namely those of dimension at most $n$, is considerably more tractable. Classical examples include the reduced and unreduced Burau representations $\Bur_{t}$ and $\UBur_{t}$, originating from Burau’s work~\cite{burau}, and the Tong–Yang–Ma representation $\TYM_{t}$ \cite{TYM}; see \S\ref{subsec:tym}, \S\ref{sec:bur}, and \S\ref{ss:reduced-Burau} for their definitions and basic properties. In particular, the classification results of Formanek, Sysoeva, and Formanek–Lee–Sysoeva–Vazirani \cite{Formanek,Formanek-Lee-Sysoeva-Vazirani,Sysoeva1,Sysoeva2} establish the central role of the Burau and Tong–Yang–Ma representations among irreducible low-dimensional representations of $\B_{n}$ with ground ring $R=\bC$. Another fundamental example is the Lawrence–Krammer–Bigelow representation $\LKB_{t,q}$; see \S\ref{ss:LKB-representation}. Arising from the homological representations introduced by Lawrence~\cite{Lawrence1}, its faithfulness was established independently by Bigelow~\cite{bigelow2001braid} and Krammer~\cite{KrammerLK}, thereby proving that braid groups are linear.

\paragraph*{Specialisations.}
The representations considered in this paper are naturally defined over group rings of free abelian groups. More precisely, their coefficient rings are of the form $R=\bZ[Q]$, where $Q\cong\bZ^{r}$, with $R=\bZ[t^{\pm1}]$ for the Burau, unreduced Burau, and Tong--Yang--Ma representations, and $R=\bZ[t^{\pm1},q^{\pm1}]$ for the Lawrence--Krammer--Bigelow representation. A choice of parameters $\bzeta=(\zeta_{1},\ldots,\zeta_{r})\in(\bC^{\times})^{r}$ determines a ring homomorphism
\[
\Xi_{\bzeta}\colon
\bZ[Q]\longrightarrow\bC.
\]
extending the canonical inclusion $\bZ\hookrightarrow\bC$ and sending $t_{i}$ to $\zeta_{i}$ for each $1\leq i\leq r$. We regard $\bC$ as a $\bZ[Q]$-module through $\Xi_{\bzeta}$. If $V$ is a representation of $\B_{n}$ over $\bZ[Q]$, we define its \emph{specialisation at} $\bzeta$ to be the complex representation
\begin{equation}\label{eq:def-specialisation}
V_{\bzeta}
:=
V\otimes_{\bZ[Q]}\bC_{\bzeta},    
\end{equation}
where $\bC_{\bzeta}:=\bC[Q]/(t_{1}-\zeta_{1},\ldots,t_{r}-\zeta_{r})\cong\bC$ is regarded as a $\bZ[Q]$-module via $\Xi_{\bzeta}$. An important observation for our computations is that taking coinvariants commutes with extension of scalars. Indeed, since coinvariants are defined as a cokernel and tensor products preserve cokernels, there is a natural isomorphism
\[
(V_{\bzeta})_{\B_{n}} \cong (V_{\B_{n}})\otimes_{\bZ[Q]}\bC_{\bzeta}.
\]
In particular, the coinvariants of the complex specialisations can be obtained directly from those of the original representations, without any flatness assumption. By contrast, taking invariants does not, in general, commute with extension of scalars.

Studying these specialisations is of particular interest, as the structure of a representation may change at certain parameter values. Properties such as irreducibility, and the dimensions of its invariant and coinvariant spaces, may hold generically but fail on exceptional subsets of the parameter space $(\bC^{\times})^{r}$. One of the purposes of the computations below is precisely to determine these exceptional parameter values for the families of representations considered here.

\paragraph*{Invariants and coinvariants.} For a representation $\rho\colon \B_{n}\to\Aut_{R}(V)$, recall that the module of $\B_{n}$-invariants is defined as
\[
V^{\B_{n}}
=
\{v\in V: \rho(\sigma)(v)=v \text{ for each }\sigma\in\B_{n}\}.
\]
Dually, the module of $\B_{n}$-coinvariants is defined as
\[
V_{\B_{n}}
=
V\big/\langle \rho(\sigma)(v)-v : \sigma\in\B_{n},\ v\in V\rangle.
\]
See for instance \cite[\S II.2]{Brown} for further details. The main purpose of this paper is to give systematic computations of the invariants and coinvariants of several fundamental families of braid group representations: the reduced and unreduced Burau, Tong--Yang--Ma, and Lawrence--Krammer--Bigelow representations, their duals as well as certain tensor products and tensor powers of these representations. We first carry out these computations over their natural coefficient rings, which are Laurent polynomial rings over $\bZ$. In particular, our arguments are integral and do not require the coefficient ring to be a field or, more generally, a principal ideal domain. We then determine the invariants and coinvariants of their complex specialisations.

\begin{notation}
To express the main results of this paper, we fix a few notations. We write $R:=\bZ[t^{\pm1},q^{\pm1}]$ and $R_{d}:=\bZ[t_{1}^{\pm1},\dots,t_{d}^{\pm1}]$.

For $n\geq2$, let
\[
\mu_{n}^{*}:=\{\zeta\in\bC^{\times}:\zeta^{n}=1,\ \zeta\neq1\}.
\]

Let $\cP_{d}$ denote the set of partitions of $[d]:=\{1,\dots,d\}$, i.e.~each element $\{p_{1},p_{2},\dots,p_{m}\}$ of $\cP_{d}$ consists of disjoint subsets $p_{i}$ of $[d]$ such that $\cup_{i} p_{i} = [d]$. For each nonempty subset $p_{i}\subseteq[d]$, we define
\[
t_{p_{i}}:=\prod_{k\in p_{i}}t_{k}.
\]

Furthermore, for $\bzeta=(\zeta_{1},\dots,\zeta_{d})\in(\bC^{\times})^{d}$, set
\[
m_{d}(\bzeta)
:=
\#\left\{
\{p_{1},\dots,p_{m}\}\in\cP_{d}:
\zeta_{p_{i}}=1
\text{ for all }1\leq i\leq m
\right\},
\]
where $\zeta_{p_{i}}=\prod_{k\in p_{i}}\zeta_{k}$.
\end{notation}

\begin{athm}\label{thm:mainA}
The main computations of this paper are summarised in the following tables. We give both invariants and coinvariants for the individual representations and their complex specialisations, and restrict to coinvariants for tensor products. 

\medskip
\noindent\textbf{Individual representations.}

\begin{center}
\begingroup
\small
\def\arraystretch{1.25}
\setlength{\tabcolsep}{3pt}
\begin{tabular}{!{\vrule width 1.5pt}c|c|c|c||c!{\vrule width 1.5pt}}
\noalign{\hrule height 1.5pt}
\textbf{Representations}
& \textbf{Invariants}
& \textbf{Coinvariants}
& \textbf{Bound}
& \textbf{Ref.}\\
\noalign{\hrule height 1.5pt}

$\TYM_{t}$, $\TYM_{t}^{\vee}$
& $0$
& $\bZ[t^{\pm1}]/(1-t)$
& $n\geq2$
& Prop.~\ref{prop:coinvariants-TYM}, Cor.~\ref{coro:invariants-TYM}\\
\noalign{\hrule height 1.5pt}

$\UBur_{t}$, $\UBur_{t}^{\vee}$
& $\bZ[t^{\pm1}]$
& $\bZ[t^{\pm1}]$
& $n\geq2$
& Prop.~\ref{prop:coinvariants-UBur}, Cor.~\ref{cor:invariants-UBur}\\
\noalign{\hrule height 1.5pt}

$\Bur_{t}$
& $0$
& $0$
& $n\geq3$
& Prop.~\ref{prop:coinvariants-Bur}, Cor.~\ref{coro:invariants-Bur}\\
\cline{1-5}

$\Bur_{t}^{\vee}$
& $0$
& $\dfrac{\bZ[t^{\pm1}]}{(1+t+\cdots+t^{n-1})}$
& $n\geq3$
& Cor.~\ref{coro:invariants-Bur-dual}, Prop.~\ref{prop:coinvariants-BurDual}\\
\noalign{\hrule height 1.5pt}

$\LKB_{t,q}$
& $0$
& $R/(1+t,1-q)$
& $n\geq3$
& Prop.~\ref{prop:coinvariants-LKB}, Cor.~\ref{coro:invariants-LKB}\\
\cline{1-5}

$\LKB_{t,q}^{\vee}$
& $0$
& $R/(1+tq^{n})$
& $n\geq2$
& Cor.~\ref{coro:invariants-LKB-dual}, Prop.~\ref{prop:coinvariants-LKB-dual}\\
\noalign{\hrule height 1.5pt}

\end{tabular}
\endgroup
\end{center}

\medskip
\noindent\textbf{Complex specialisations.}

\begin{center}
\begingroup
\small
\def\arraystretch{1.25}
\setlength{\tabcolsep}{3pt}
\begin{tabular}{!{\vrule width 1.5pt}c|c|c|c||c!{\vrule width 1.5pt}}
\noalign{\hrule height 1.5pt}
\textbf{Representations}
& \textbf{Invariants}
& \textbf{Coinvariants}
& \textbf{Bound}
& \textbf{Ref.}\\
\noalign{\hrule height 1.5pt}

$\TYM_{\zeta}$, $\TYM_{\zeta}^{\vee}$
&
$\begin{cases}
\bC & \zeta=1,\\
0 & \text{otherwise}
\end{cases}$
&
$\begin{cases}
\bC & \zeta=1,\\
0 & \text{otherwise}
\end{cases}$
& $n\geq2$
& Prop.~\ref{prop:coinvariants-TYM}, Cor.~\ref{coro:invariants-TYM}\\
\noalign{\hrule height 1.5pt}

$\UBur_{\zeta}$, $\UBur_{\zeta}^{\vee}$
& $\bC$
& $\bC$
& $n\geq2$
& Prop.~\ref{prop:coinvariants-UBur}, Cor.~\ref{cor:invariants-UBur}\\
\noalign{\hrule height 1.5pt}

$\Bur_{\zeta}$
&
$\begin{cases}
\bC & \zeta\in\mu_{n}^{*},\\
0 & \text{otherwise}
\end{cases}$
& $0$
& $n\geq3$
& Prop.~\ref{prop:coinvariants-Bur}, Cor.~\ref{coro:invariants-Bur}\\
\cline{1-5}

$\Bur_{\zeta}^{\vee}$
& $0$
&
$\begin{cases}
\bC & \zeta\in\mu_{n}^{*},\\
0 & \text{otherwise}
\end{cases}$
& $n\geq3$
& Cor.~\ref{coro:invariants-Bur-dual}, Prop.~\ref{prop:coinvariants-BurDual}\\
\noalign{\hrule height 1.5pt}

$\LKB_{\bzeta}$
&
$\begin{cases}
\bC & \zeta_{t}\zeta_{q}^{n}=-1,\\
0 & \text{otherwise}
\end{cases}$
&
$\begin{cases}
\bC & \zeta_{t}=-1,\ \zeta_{q}=1,\\
0 & \text{otherwise}
\end{cases}$
& $n\geq3$
& Prop.~\ref{prop:coinvariants-LKB}, Cor.~\ref{coro:invariants-LKB}\\
\cline{1-5}

$\LKB_{\bzeta}^{\vee}$
&
$\begin{cases}
\bC & \zeta_{t}=-1,\ \zeta_{q}=1,\\
0 & \text{otherwise}
\end{cases}$
&
$\begin{cases}
\bC & \zeta_{t}\zeta_{q}^{n}=-1,\\
0 & \text{otherwise}
\end{cases}$
& $n\geq3$
& Cor.~\ref{coro:invariants-LKB-dual}, Prop.~\ref{prop:coinvariants-LKB-dual}\\
\noalign{\hrule height 1.5pt}

\end{tabular}
\endgroup
\end{center}

\medskip
\noindent\textbf{Tensor products.}

\begin{center}
\begingroup
\small
\def\arraystretch{1.35}
\setlength{\tabcolsep}{3pt}
\begin{tabular}{!{\vrule width 1.5pt}c|c|c|c||c!{\vrule width 1.5pt}}
\noalign{\hrule height 1.5pt}
\textbf{Representations}
& \shortstack{\textbf{Over canonical coefficient ring}}
& \textbf{Complex specialisations}
& \textbf{Bound}
& \textbf{Ref.}\\
\noalign{\hrule height 1.5pt}

$\TYM_{t_{1}}\otimes\cdots\otimes\TYM_{t_{d}}$
&
$\displaystyle
\bigoplus_{\{p_{1},\dots,p_{m}\}\in\cP_{d}}
\frac{R_{d}}{(1-t_{p_{1}},\dots,1-t_{p_{m}})}$
&
$\bC^{m_{d}(\bzeta)}$
& $n\geq d+1$
& Th.~\ref{thm:coinvariants-TYM-tensor}\\
\noalign{\hrule height 1.5pt}

$\UBur_{t}\otimes\UBur_{q}$
& $R\oplus R/(1-tq)$
&
$\begin{cases}
\bC^{2} & \zeta_{t}\zeta_{q}=1,\\
\bC & \text{otherwise}
\end{cases}$
& $n\geq2$
& Prop.~\ref{prop:coinvariants-UBur2}\\
\noalign{\hrule height 1.5pt}

$\TYM_{t}\otimes\UBur_{q}$
&
$\dfrac{R}{(1-t)}
\oplus\dfrac{R}{(1-t,1-q)}$
&
$\begin{cases}
\bC^{2} & \zeta_{t}=\zeta_{q}=1,\\
\bC & \zeta_{t}=1,\ \zeta_{q}\neq1,\\
0 & \zeta_{t}\neq1
\end{cases}$
& $n\geq3$
& Prop.~\ref{prop:coinvariants-TYMxUBur}\\
\noalign{\hrule height 1.5pt}

$\Bur_{t}\otimes\Bur_{q}$
& $R/(1-tq)$
&
$\begin{cases}
\bC & \zeta_{t}\zeta_{q}=1,\\
0 & \text{otherwise}
\end{cases}$
& $n\geq3$
& Prop.~\ref{prop:coinvariants-Bur2}\\
\noalign{\hrule height 1.5pt}

$\TYM_{t}\otimes\Bur_{q}$
& $R/(1-t,1-q)$
&
$\begin{cases}
\bC & \zeta_{t}=\zeta_{q}=1,\\
0 & \text{otherwise}
\end{cases}$
& $n\geq3$
& Prop.~\ref{prop:coinvariants-TYMxBur}\\
\noalign{\hrule height 1.5pt}

$\UBur_{t}\otimes\Bur_{q}$
& $R/(1-tq)$
&
$\begin{cases}
\bC & \zeta_{t}\zeta_{q}=1,\\
0 & \text{otherwise}
\end{cases}$
& $n\geq3$
& Prop.~\ref{prop:coinvariants-UBurxBur}\\
\noalign{\hrule height 1.5pt}
\end{tabular}
\endgroup
\end{center}
\end{athm}

Note that the invariants of tensor products can be investigated using Lemma~\ref{lem:invariants-coinvariants-duality} and computations for their dual representations; this is work in progress.

Some of the computations in Theorem~\ref{thm:mainA} recover earlier results on the (co)homology of braid groups with twisted coefficients, obtained by different methods. For example:
\begin{itemizeb}
    \item \textbf{The reduced Burau representation.}
    Chen~\cite{Chen} computed the homology of braid groups with coefficients in the reduced Burau representation over $\bC[t^{\pm1}]$, as well as its complex specialisations. In particular, the degree-zero case of~\cite[Th.~1]{Chen} gives the vanishing of the coinvariants of $\Bur_{t}$ for $n\geq 3$, while for $n=2$ they are isomorphic to $\bC[t^{\pm1}]/(1+t)$. The computations for complex specialisations agree with~\cite[Prop.~7]{Chen}.

    \item \textbf{The Tong--Yang--Ma representation.}
    Callegaro, Moroni, and Salvetti~\cite{CallegaroMoroniSalvetti} computed the cohomology of braid groups with coefficients in the Tong--Yang--Ma representation over $\bQ[t^{\pm1}]$. In particular, their results imply that, for $n\geq 3$, the coinvariants of $\TYM_{t}$ are isomorphic to $\bQ[t^{\pm1}]/(1-t)$. The authors point out that their computations could be partially extended to integral coefficients, but that substantial difficulties arise because $\bZ[t^{\pm1}]$ is not a principal ideal domain. Our direct computation, by contrast, applies over $\bZ[t^{\pm1}]$ and establishes the corresponding result integrally, for each $n\geq 2$.
\end{itemizeb}
More generally, Theorem~\ref{thm:mainA} provides a systematic computation of the invariants and coinvariants of the representations considered in this paper, including their duals and complex specialisations. In particular, we determine precisely the values of the parameters for which the invariants or coinvariants are nonzero, thereby identifying the exceptional specialisations at which their dimensions jump. This gives a complete description of the degree-zero behaviour of these families, including cases that are not covered by the preceding results.

\paragraph*{Reducibility results.}
Our computations also have consequences for the reducibility of the representations considered in this paper. Indeed, let $R$ be a commutative ring and let $V$ be an $R[\B_{n}]$-module. The module $V^{\B_{n}}$ is the largest submodule of $V$ on which $\B_{n}$ acts trivially, whereas $V_{\B_{n}}$ is the largest quotient of $V$ on which $\B_{n}$ acts trivially. Thus, invariants give rise to subrepresentations with trivial action, while coinvariants give rise to quotient representations with trivial action. In particular, when $V$ is free of rank greater than one, either condition implies that $V$ is reducible. Theorem~\ref{thm:mainA} therefore allows us to determine explicitly these trivial subrepresentations and quotients, both over the natural Laurent polynomial coefficient rings and after complex specialisation. In the latter case, our computations provide reducibility criteria in terms of the specialisation parameters, including exceptional specialisations where the dimensions of the invariant or coinvariant spaces exceed one.

\begin{coro}
For the representations and bounds considered in Theorem~\ref{thm:mainA}, the following hold.

\begin{itemizeb}
    \item Over their natural coefficient rings, the unreduced Burau representations and their duals, together with $\UBur_{t}\otimes\UBur_{q}$, have nonzero subrepresentations with trivial action. All the other representations listed in the theorem have no nonzero subrepresentation with trivial action.

    \item All these representations, except $\Bur_{t}$, admit a nonzero quotient with trivial action.

    \item After complex specialisation, subrepresentations and quotients with trivial action occur precisely at the parameter values for which the corresponding invariant and coinvariant spaces in Theorem~\ref{thm:mainA} are nonzero.
\end{itemizeb}
In particular, each such nonzero invariant or coinvariant space yields a reducibility criterion whenever the representation has dimension greater than one.
\end{coro}

\begin{rmk}
The reducibility of several of the complex specialisations of these representations has already been extensively studied. In particular, Formanek~\cite[Lem.~6]{Formanek} determined the irreducibility criterion for complex specialisations of the reduced Burau representation; see also~\cite[Lem.~6.1]{Formanek-Lee-Sysoeva-Vazirani}. The corresponding results for the Tong--Yang--Ma representation follow from the classification of low-dimensional irreducible representations of braid groups~\cite{Sysoeva1,Formanek-Lee-Sysoeva-Vazirani,Sysoeva2}. For the complex specialisations of the Lawrence--Krammer--Bigelow representation, generic irreducibility was established by Zinno~\cite{Zinno}, while reducibility at exceptional complex parameters was investigated by Levaillant and Wales~\cite{Levaillant-Irred,LevaillantWales}. The corresponding invariant subspaces were subsequently classified by Levaillant~\cite{Levaillant-invariants}, under suitable assumptions on the parameters.

Our approach is complementary to these representation-theoretic results. Rather than classifying all invariant subspaces, we determine the maximal subrepresentations and quotients with trivial action over the natural coefficient rings, as well as their dimensions and the precise parameter values at which they arise after complex specialisation. In particular, our computations for $\LKB_{t,q}$ and its dual are carried out over $\bZ[t^{\pm1},q^{\pm1}]$, and those for their complex specialisations apply to all nonzero values of the parameters, without the semisimplicity assumptions imposed in the aforementioned classification results.
\end{rmk}

\begin{rmk}\label{rmk:non-splitting}
The existence of a subrepresentation with trivial action does not necessarily imply that it splits off as a direct summand, even when both the subrepresentation and the corresponding quotient with trivial action are free over the coefficient ring. For instance, the invariant submodule of $\UBur_{t}$ is freely generated over $R=\bZ[t^{\pm1}]$ by $u_{t}=v_{1}+tv_{2}+\cdots+t^{n-1}v_{n}$, while its quotient with trivial action is induced by the $R$-linear map $\varepsilon_{t}$ defined by $\varepsilon_{t}(v_{j})=1$ for each $1\leq j\leq n$. Since $\varepsilon_{t}(u_{t})=1+t+\cdots+t^{n-1}$, is not invertible in $R$, the subrepresentation with trivial action does not split off over $R$. After complex specialisation at $\zeta\in\bC^{\times}$, the subrepresentation with trivial action of $\UBur_{\zeta}$ splits off if and only if $1+\zeta+\cdots+\zeta^{n-1}\neq0$. In particular, it does not split when $\zeta\neq1$ is an $n^\mathrm{th}$ root of unity.

Similarly, let $\bzeta=(\tau,\xi)\in(\bC^{\times})^{2}$ satisfy $\tau\xi^{n}=-1$ and $\xi\neq1$. Then $\LKB_{\bzeta}$ contains a one-dimensional subrepresentation with trivial action but admits no nonzero quotient with trivial action. Therefore, this subrepresentation cannot split off as a direct summand.
\end{rmk}

\paragraph*{Homological stability.}
Our results may also be viewed from the perspective of twisted homological stability. For a representation $\rho\colon \B_{n}\to\Aut_{R}(V)$, recall that the module of invariants is the zeroth cohomology group $V^{\B_{n}}=H^{0}(\B_{n};V)$, whereas the module of coinvariants is the zeroth homology group $V_{\B_{n}}=H_{0}(\B_{n};V)$. See for instance \cite[\S II.2]{Brown} for further details. Hence, our computations of Theorem~\ref{thm:mainA} determine the twisted homology and cohomology of the braid groups in degree zero for the coefficients given by the families of representations considered in this paper. Randal-Williams and Wahl \cite[Th.~D, Th.~5.22]{RWW} prove homological stability for braid groups with polynomial coefficient systems. Namely, if $F$ is a coefficient system of degree $r$ at $N$, their result implies that the stabilisation morphism
\begin{equation}\label{eq:stabilisation-map}
\St_{n,i}\colon H_{i}(\B_{n};F(n)) \longrightarrow H_{i}(\B_{n+1};F(n+1)).
\end{equation}
is an isomorphism for $i\leq \frac{n-2}{2}-r$, provided that $n>N$. Several of the families of braid group representations considered in this paper define polynomial coefficient systems to which this result applies. More precisely, the unreduced Burau and Tong--Yang--Ma representations $\UBur_{t}$ and $\TYM_{t}$, together with their complex specialisations, define coefficient systems of degree $1$ at $0$ by~\cite[Proposition~3.25]{S1}; see also~\cite[Example~4.15]{RWW} for the earlier result concerning $\UBur_{t}$. Similarly, the reduced Burau representation $\Bur_{t}$ and its complex specialisations define coefficient systems of degree $1$ at $1$ by~\cite[Proposition~3.28]{S1}, while the Lawrence--Krammer--Bigelow representation $\LKB_{t,q}$ and its complex specialisations define coefficient systems of degree $2$ at $0$ by~\cite[Proposition~3.33]{S1}. Moreover, the $p$-fold tensor power of a coefficient system of degree $d$ is a coefficient system of degree at most $pd$; see, for instance,~\cite[Corollary~1.5]{AnnexLS1}.

Since we compute their zeroth homology explicitly, we can compare the resulting stability behaviour with the general bounds of Randal-Williams and Wahl~\cite{RWW}. In degree zero, their theorem guarantees that the stabilisation map is an isomorphism for $n\geq 4$ for coefficient systems of degree $1$ (at $0$ or $1$). This applies to $\TYM_{t}$, $\UBur_{t}$, and $\Bur_{t}$, together with their complex specialisations. For coefficient systems of degree $2$ at $0$, such as $\LKB_{t,q}$ and tensor products of two degree-$1$ coefficient systems, the corresponding bounds for isomorphisms are $n\geq 6$. More generally, a tensor product of $d$ of the above degree-$1$ coefficient systems (namely, reduced or unreduced Burau or Tong--Yang--Ma coefficient systems, including their complex specialisations), so the general stability theorem guarantees an isomorphism in degree zero for $n\geq 2d+2$. Our explicit computations allow us to investigate whether these general bounds are optimal for the representations considered here.

\begin{coro}\label{coro:bound-homological-stability}
For each family of representations $V=\{V_{n}\}_{n\geq2}$ listed below, the stabilisation morphism $\St_{n,0}$ is an isomorphism for all $n\geq N$, where $N$ is given in the second column. The third column gives the corresponding bound from \cite[Theorems~D and~5.22]{RWW}.
\[
\begin{array}{l|c|c}
\text{Coefficient system}
& N
& \text{Bound from \cite{RWW}}\\
\hline
\TYM_{t} & 2 & 4\\
\UBur_{t} & 2 & 4\\
\Bur_{t} & 3 & 4\\
\LKB_{t,q} & 3 & 6\\
\TYM_{t_{1}}\otimes\cdots\otimes\TYM_{t_{d}}
& d+1 & 2d+2\\
\UBur_{t}\otimes\UBur_{q} & 2 & 6\\
\TYM_{t}\otimes\UBur_{q} & 3 & 6\\
\Bur_{t}\otimes\Bur_{q} & 3 & 6\\
\TYM_{t}\otimes\Bur_{q} & 3 & 6\\
\UBur_{t}\otimes\Bur_{q} & 3 & 6
\end{array}
\]
The same bounds hold for the corresponding complex specialisations. The bounds for $\TYM_{t}$, $\UBur_{t}$, $\Bur_{t}$, and $\LKB_{t,q}$ are sharp.
\end{coro}

In particular, our explicit computations improve the general bounds of~\cite{RWW} for all these families.

There is also a complementary phenomenon. Some of the families studied below do not satisfy homological stability, even in degree zero: their coinvariants depend essentially on $n$ and therefore cannot eventually be identified by the stabilisation morphisms. Our explicit computations thus provide obstructions to twisted homological stability which are already visible at the level of $H_{0}$.

\begin{coro}\label{coro:non-homological-stability}
The families of dual reduced Burau representations $\Bur_{t}^{\vee}$ and dual Lawrence--Krammer--Bigelow representations $\LKB_{t,q}^{\vee}$ do not satisfy homological stability in degree zero.
\end{coro}

\begin{rmk}
Corollary~\ref{coro:non-homological-stability} is particularly interesting in view of the homological stability results of Randal-Williams and Wahl~\cite{RWW}. Indeed, the ranks of $\Bur_{t}^{\vee}$ and $\LKB_{t,q}^{\vee}$ grow polynomially with $n$, being $n-1$ and $n(n-1)/2$, respectively. Nevertheless, neither family can be assembled into a coefficient system of finite degree in the sense of~\cite{RWW}, since the existence of such a coefficient system would imply eventual homological stability in degree zero. In particular, polynomial growth of the ranks of a family of representations does not imply polynomiality as a coefficient system. We suspect that the obstruction already arises at the level of functoriality: these families may fail to admit the structure of coefficient systems in the sense of~\cite{RWW}, even before considering the question of finite degree or polynomiality.
\end{rmk}

\paragraph*{Toolbox.} The computations in this paper rely on two elementary descriptions of invariants and coinvariants. Since the standard generators $\sigma_{1},\ldots,\sigma_{n-1}$ generate $\B_{n}$, for any representation $\rho\colon\B_{n}\to\Aut_{R}(V)$, we define the $R$-linear maps
\[
\begin{aligned}
\Phi_{\rho}\colon V
&\longrightarrow \bigoplus_{i=1}^{n-1}V,\\
v
&\longmapsto
\bigl(v-\rho(\sigma_{1})(v),\ldots,
v-\rho(\sigma_{n-1})(v)\bigr),\\[1ex]
\Psi_{\rho}\colon \bigoplus_{i=1}^{n-1}V
&\longrightarrow V,\\
(v_{1},\ldots,v_{n-1})
&\longmapsto
\sum_{i=1}^{n-1}
\bigl(v_{i}-\rho(\sigma_{i})(v_{i})\bigr).
\end{aligned}
\]
The invariants and coinvariants of $V$ are then given by
\begin{equation}\label{eq:def-invariants-coinvariants-Phi-Psi}
V^{\B_{n}}=\ker(\Phi_{\rho})
\qquad\text{and}\qquad
V_{\B_{n}}=V/\Image(\Psi_{\rho})
=\operatorname{coker}(\Psi_{\rho}).
\end{equation}
Thus, the coinvariants are obtained by imposing the relations $\rho(\sigma_{i})(v)=v$, whereas the invariants are obtained by solving these relations simultaneously. In particular, both may be computed directly from the matrices of the standard braid generators.

Moreover, taking the inverse of a representation does not change its invariants or coinvariants, since the elements of $\B_{n}$ and their inverses impose the same relations. For a braid representation $\rho\colon \B_{n} \to \Aut_{R}(V)$, let $V^{\tr}$, $V^{\inv}$, and $V^{\vee}$ denote the underlying $\B_{n}$-modules for the transpose, inverse, and dual representations of $\rho$, respectively. It follows that, as $V^{\vee}=(V^{\tr})^{\inv}$, we have natural isomorphisms
\begin{equation}\label{eq:transpose-dual-invariants-coinvariants}
V_{\B_{n}}^{\inv}\cong V_{\B_{n}},
\qquad
(V^{\inv})^{\B_{n}}\cong V^{\B_{n}},
\qquad
(V^{\tr})_{\B_{n}}\cong (V^{\vee})_{\B_{n}},
\qquad
(V^{\tr})^{\B_{n}}\cong (V^{\vee})^{\B_{n}}.
\end{equation}
In particular, the computations for the transpose representation also determine the invariants and coinvariants of the dual representation, and conversely.

Furthermore, a useful tool for our computations is the duality between invariants and coinvariants, recalled in Lemma~\ref{lem:invariants-coinvariants-duality}. Since all the representations considered in this paper are free of finite rank over their respective coefficient rings, this duality allows us to recover their invariants from the coinvariants of their duals. Our general strategy will therefore be to compute first the coinvariants $V_{\B_{n}}$ and $(V^{\vee})_{\B_{n}}$, and then deduce the invariants $(V^{\vee})^{\B_{n}}$ and $V^{\B_{n}}$, respectively.

%\paragraph*{Perspectives.} As explained above, the invariants and coinvariants of a $\B_{n}$-representation $V$ can be interpreted as the degree-zero cohomology and homology groups $H^{0}(\B_{n};V)$ and $H_{0}(\B_{n};V)$, respectively. The present work thus provides the starting point for the study of braid group (co)homology with coefficients in these representations. In the forthcoming paper~\cite{hoangsoulie1}, we extend these computations to degree one, determining $H^{1}(\B_{n};V)$ and $H_{1}(\B_{n};V)$ using the Fox--Neuwirth /Fuks chain and cochain complexes.

\paragraph*{Acknowledgements.}
We thank Nir Gadish, Ben Knudsen, Martin Palmer, and Craig Westerland for helpful discussions. We are particularly grateful to Calista Bernard for bringing us together and initiating this collaboration. A. S. acknowledges support from the French National Research Agency (ANR) through the project HoReStHo (ANR-26-CE40-3880).

%-------------------------------------------------------

\section{The low-dimensional representations}\label{sec:low_dim}

We begin by computing the invariants and coinvariants of the Tong--Yang--Ma, reduced, and unreduced Burau representations, together with their complex specialisations. As explained in the introduction, these are low-dimensional representations of $\B_{n}$, in the sense that their dimensions are at most $n$.

\begin{notation}
For $m\in \bZ_{>0}$, let $\{v_{1},\ldots,v_{m}\}$ be the standard basis for $V=\bZ[t^{\pm1}]^{\oplus m}$. That is,
\[
V=\bZ[t^{\pm1}]\{v_{1},\ldots,v_{m}\},
\]
where $v_{j}$ has $1$ in the $j$-th position and $0$ elsewhere.
\end{notation}

%-------------------------------------------------------

\subsection{The Tong--Yang--Ma representation}\label{subsec:tym}

%-------------------------------------------------------

Let $\TYM_{t}$ denote the Tong–Yang–Ma representation of $\B_{n}$ over $\bZ[t^{\pm1}]$. Identifying $\TYM_{t}$ with $\bZ[t^{\pm1}]^{\oplus n}$, the action of each standard generator $\sigma_{i}$ is given by
\[
\TYM_{t}(\sigma_{i})(v_{j})=
\begin{cases}
v_{j} & \text{if }j\notin\{i,i+1\},\\
v_{i} & \text{if }j=i+1,\\
tv_{i+1} & \text{if }j=i.
\end{cases}
\]
% The action of a generator $\sigma_{i}$ of $\B_{n}$ on $\TYM_{t}\cong \bZ[t^{\pm1}]^{\oplus n}$ is given by
% \[\sigma_{i} \mapsto
% \begin{pNiceArray}{c|cc|c}
%   I_{i-1} & & & \\
%   \hline
%   & 0 & 1 & \\
%   & t & 0 & \\
%   \hline
%   & & & I_{n-i-1} \\
% \end{pNiceArray}.
% \]
It is known that $\TYM_{t}\cong\TYM_{t}^{\tr}$, with an intertwining isomorphism given by conjugation either by the anti-diagonal matrix whose nonzero entries are all equal to $1$, or by the diagonal matrix $\mathrm{diag}(1,t,\ldots,t^{n-1})$. Moreover, a direct computation of the inverse matrices yields $\TYM_{t}^{\inv}\cong\TYM_{t^{-1}}^{\tr}$.
%\[
%\begin{pNiceArray}{c|cc|c}
%  I_{i-1} & & & \\
%  \hline
%  & 0 & 1 & \\
%  & t & 0 & \\
%  \hline
%  & & & I_{n-i-1} \\
%\end{pNiceArray}^{-1}
%=
%\begin{pNiceArray}{c|cc|c}
%  I_{i-1} & & & \\
%  \hline
%  & 0 & t^{-1} & \\
%  & 1 & 0 & \\
%  \hline
%  & & & I_{n-i-1} \\
%\end{pNiceArray}
%=
%\begin{pNiceArray}{c|cc|c}
%  I_{i-1} & & & \\
%  \hline
%  & 0 & 1 & \\
%  & t^{-1} & 0 & \\
%  \hline
%  & & & I_{n-i-1} \\
%\end{pNiceArray}^\mathrm{T},
%\]
%More concretely, the action of $\B_{n}$ on the inverse representation $\TYM_{t}^{\inv}$ is given by
%\[\sigma_{i}^{-1}(v_{j}) = \begin{cases} v_{j} & i \neq j-1, j\\ t^{-1}v_{j-1} & i = j-1\\ v_{j+1} & i=j. \end{cases}\]
As a result, we have
\begin{equation}\label{eq:dual-TYM}
    \TYM_{t}^{\vee} \cong \left(\TYM_{t}^{\inv}\right)^{\tr} \cong \TYM_{t^{-1}}.
\end{equation}
For $\zeta\in\bC^{\times}$, recall that $\TYM_{\zeta}$ denotes the specialisation of $\TYM_{t}$ at $t=\zeta$, as defined in~\eqref{eq:def-specialisation}. We begin by computing the coinvariants of $\TYM_{t}$, $\TYM_{t}^{\vee}$, and their complex specialisations.

\begin{prop}\label{prop:coinvariants-TYM}
For $n\geq 2$, we have
\[
\left(\TYM_{t}\right)_{\B_{n}} \cong \bigl(\TYM_{t}^{\vee}\bigr)_{\B_{n}}
\cong \bZ[t^{\pm1}]/(1-t).
\]
Moreover, for $\zeta\in\bC^{\times}$,
\[
(\TYM_{\zeta})_{\B_{n}}
\cong
\bigl(\TYM_{\zeta}^{\vee}\bigr)_{\B_{n}}
\cong
\begin{cases}
0 & \text{if }\zeta\neq 1,\\
\bC & \text{if }\zeta=1.
\end{cases}
\]
\end{prop}

\begin{proof}
We compute from the definition of the representation $\TYM_{t}$ that 
\[
(1-\TYM_{t}(\sigma_{i}))(v_{j})=
\begin{cases}
0 & \text{if }i\neq j-1,j,\\
v_{j}-v_{j-1} & \text{if }i=j-1,\\
v_{j}-tv_{j+1} & \text{if }i=j.
\end{cases}
\]
Hence, we deduce that
\[
\Image(\Psi_{\TYM_{t}})
=
\bZ[t^{\pm1}]
\bigl\{
v_{j+1}-v_{j},\,
v_{j}-tv_{j+1}
: 1\leq j\leq n-1
\bigr\}.
\]
Therefore, using the description of coinvariants in
\eqref{eq:def-invariants-coinvariants-Phi-Psi}, we obtain
\[
\begin{aligned}
(\TYM_{t})_{\B_{n}}
&\cong
\frac{\bZ[t^{\pm1}]\{v_{1},\ldots,v_{n}\}}
{\bigl(v_{j+1}-v_{j},\;v_{j}-tv_{j+1}:1\leq j\leq n-1\bigr)}
\\
&\cong
\frac{\bZ[t^{\pm1}]\{v_{1},\ldots,v_{n}\}}
{\bigl(v_{j+1}-v_{1},\;(1-t)v_{1}:1\leq j\leq n-1\bigr)}
\\
&\cong
\frac{\bZ[t^{\pm1}]\{v_{1}\}}{((1-t)v_{1})}
\cong
\bZ[t^{\pm1}]/(1-t).
\end{aligned}
\]

We now consider the complex specialisation $\TYM_{\zeta}$. If $\zeta\neq 1$, the same argument as above gives $(\TYM_{\zeta})_{\B_{n}}\cong\bC/(1-\zeta)=0$. If $\zeta=1$, we have
\[
\Image(\Psi_{\TYM_{1}})
=
\bC\{v_{j}-v_{j+1}:1\leq j\leq n-1\}.
\]
Hence, by~\eqref{eq:def-invariants-coinvariants-Phi-Psi}, we obtain
\[
(\TYM_{1})_{\B_{n}}
\cong
\frac{\bC\{v_{1},\ldots,v_{n}\}}
{(v_{j}-v_{j+1}:1\leq j\leq n-1)}
\cong
\frac{\bC\{v_{1},\ldots,v_{n}\}}
{(v_{1}-v_{j+1}:1\leq j\leq n-1)}
\cong \bC\{v_{1}\}
\cong \bC.
\]
Finally, by the isomorphisms~\eqref{eq:dual-TYM}, the coinvariants of $\TYM_{t}^{\vee}$ are obtained from those of $\TYM_{t}$ by replacing $t$ with $t^{-1}$. Since $1-t^{-1}=-t^{-1}(1-t)$ and $t$ is invertible, the ideals $(1-t^{-1})$ and $(1-t)$ coincide, yielding the desired isomorphism.
\end{proof}

We next compute the invariants of these representations. To this end, we use the following general lemma, which allows us to recover the invariants of a representation from the coinvariants of its dual.

\begin{lem}\label{lem:invariants-coinvariants-duality}
Let $R$ be a commutative ring, let $G$ be a group, and let $V$ be an $R[G]$-module which is finitely generated and projective as an $R$-module. Then there is a natural isomorphism
\begin{equation}\label{eq:invariants-dual-coinvariants}
V^{G} \cong \Hom_{R}\bigl((V^{\vee})_{G},R\bigr).
\end{equation}
\end{lem}

\begin{proof}
We first observe that an $R$-linear map $\varphi\colon V\to R$ is $G$-invariant if and only if $\varphi(gv-v)=0$ for each $g\in G$ and $v\in V$, which is equivalent to saying that $\varphi$ factors uniquely through the quotient $V\to V_{G}$. This provides a natural isomorphism
\[
(V^{\vee})^{G}
\cong
\Hom_{R}(V_{G},R).
\]
Applying this observation to $V^{\vee}$ gives $(V^{\vee\vee})^{G} \cong \Hom_{R}\bigl((V^{\vee})_{G},R\bigr)$. But since $V$ is finitely generated and projective over $R$, the canonical map $V\to V^{\vee\vee}$ is an isomorphism of $R[G]$-modules, and the result follows.
\end{proof}

The following result is then obtained by combining Lemma~\ref{lem:invariants-coinvariants-duality}, the isomorphisms~\eqref{eq:dual-TYM}, and Proposition~\ref{prop:coinvariants-TYM}:
\begin{coro}\label{coro:invariants-TYM}
For $n\geq 2$, we have $\TYM_{t}^{\B_{n}} \cong \bigl(\TYM_{t}^{\vee}\bigr)^{\B_{n}} =0$. Moreover, for $\zeta\in\bC^{\times}$, we have
\[
\TYM_{\zeta}^{\B_{n}}
\cong
\bigl(\TYM_{\zeta}^{\vee}\bigr)^{\B_{n}}
\cong
\begin{cases}
0 & \text{if }\zeta\neq 1,\\
\bC & \text{if }\zeta=1.
\end{cases}
\]
\end{coro}

%-------------------------------------------------------

\subsection{The unreduced Burau representation}\label{sec:bur}

Let $\UBur_{t}$ denote the unreduced Burau representation of $\B_{n}$ over $\bZ[t^{\pm1}]$. Identifying $\UBur_{t}$ with $\bZ[t^{\pm1}]^{\oplus n}$, the action of each standard generator $\sigma_{i}$ is given by
\[
\UBur_{t}(\sigma_{i})(v_{j})=
\begin{cases}
v_{j} & \text{if }j\notin\{i,i+1\},\\
v_{i} & \text{if }j=i+1,\\
(1-t)v_{i}+tv_{i+1} & \text{if }j=i.
\end{cases}
\]
% Alternatively, in matrix form, we have
% \[\sigma_{i} \mapsto
% \begin{pNiceArray}{c|cc|c}
%   I_{i-1} & & & \\
%   \hline
%   & 1-t & 1 & \\
%   & t & 0 & \\
%   \hline
%   & & & I_{n-i-1} \\
% \end{pNiceArray}.
% \]
Moreover, $\UBur_{t}\cong\UBur_{t}^{\tr}$, with an intertwining isomorphism given by conjugation by the diagonal matrix $\mathrm{diag}(1,t,\ldots,t^{n-1})$. Therefore, there is an isomorphism 
\begin{equation}\label{eq:UBur-dual-iso-inv}
    \UBur_{t}^{\vee}\cong\UBur_{t}^{\inv}.
\end{equation}
%Explicitly, the action of $\B_{n}$ on the inverse $\UBur_{t}^{\inv}$ is given by
%\[\sigma_{i}^{-1}(v_{j}) = \begin{cases}
%    v_{j} & i \neq j-1, j\\
%    t^{-1}v_{j-1}+(1-t^{-1})v_{j} & i = j-1\\
%    v_{j+1} & i=j.
%\end{cases}\]
For $\zeta\in\bC^{\times}$, recall that $\UBur_{\zeta}$ denotes the specialisation of $\UBur_{t}$ at $t=\zeta$, as defined in~\eqref{eq:def-specialisation}.
We begin by computing the coinvariants of the unreduced Burau representation and its complex specialisations.

\begin{prop}\label{prop:coinvariants-UBur}
For $n\geq 2$, we have
\[
(\UBur_{t})_{\B_{n}}
\cong \bigl(\UBur_{t}^{\vee}\bigr)_{\B_{n}}
\cong \bZ[t^{\pm1}]
\]
and, for any $\zeta\in\bC^{\times}$,
\[
(\UBur_{\zeta})_{\B_{n}}
\cong \bigl(\UBur_{\zeta}^{\vee}\bigr)_{\B_{n}}
\cong \bC.
\]
\end{prop}
\begin{proof}
By the definition of the unreduced Burau representation, we have
\[
(1-\UBur_{t}(\sigma_{i}))(v_{j})=
\begin{cases}
0 & \text{if }i\neq j-1,j,\\
v_{j}-v_{j-1} & \text{if }i=j-1,\\
t(v_{j}-v_{j+1}) & \text{if }i=j.
\end{cases}
\]
We deduce that
\[
\begin{aligned}
\Image(\Psi_{\UBur_{t}})
&=\bZ[t^{\pm1}]
\{v_{j+1}-v_{j},\;t(v_{j}-v_{j+1}):1\leq j\leq n-1\}\\
&=\bZ[t^{\pm1}]
\{v_{j+1}-v_{j}:1\leq j\leq n-1\},
\end{aligned}
\]
where the second equality follows from the invertibility of $t$. Hence, by~\eqref{eq:def-invariants-coinvariants-Phi-Psi},
\[
(\UBur_{t})_{\B_{n}}
\cong
\frac{\bZ[t^{\pm1}]\{v_{1},\ldots,v_{n}\}}
{(v_{j+1}-v_{j}:1\leq j\leq n-1)}\\
\cong \bZ[t^{\pm1}]\{v_{1}\}
\cong \bZ[t^{\pm1}].
\]
Since taking coinvariants commutes with extension of scalars and the preceding computation is independent of $t$, it applies equally to every complex specialisation $t=\zeta\in\bC^{\times}$, yielding $(\UBur_{\zeta})_{\B_{n}}\cong\bC$.

Finally, the results for $\bigl(\UBur_{t}^{\vee}\bigr)_{\B_{n}}$ and $\bigl(\UBur_{\zeta}^{\vee}\bigr)_{\B_{n}}$ follow from the preceding computations by combining the isomorphism~\eqref{eq:UBur-dual-iso-inv} with the general relations~\eqref{eq:transpose-dual-invariants-coinvariants}.
\end{proof}

Combining Lemma~\ref{lem:invariants-coinvariants-duality}, the isomorphism~\eqref{eq:UBur-dual-iso-inv}, the general relation~\eqref{eq:transpose-dual-invariants-coinvariants}, and Proposition~\ref{prop:coinvariants-UBur}, we obtain the following result.

\begin{coro}\label{cor:invariants-UBur}
For $n\geq 2$, we have $\UBur_{t}^{\B_{n}}
\cong \bigl(\UBur_{t}^{\vee}\bigr)^{\B_{n}}
\cong \bZ[t^{\pm1}]$. Moreover, for any $\zeta\in\bC^{\times}$,
\[
\UBur_{\zeta}^{\B_{n}}
\cong \bigl(\UBur_{\zeta}^{\vee}\bigr)^{\B_{n}}
\cong \bC.
\]
\end{coro}

%-------------------------------------------------------

\subsection{The reduced Burau representation}\label{ss:reduced-Burau}

Let $\Bur_{t}$ denote the reduced Burau representation of $\B_{n}$ over $\bZ[t^{\pm1}]$. Identifying $\Bur_{t}$ with $\bZ[t^{\pm1}]^{\oplus n-1}$, the action of each standard generator $\sigma_{i}$ is given by
\[
\Bur_{t}(\sigma_{i})(v_{j}) = \begin{cases}
    v_{j} & \text{if }i \neq j-1, j, j+1,\\
    v_{j-1}+v_{j} & \text{if }i = j-1,\\
    -tv_{j} & \text{if }i=j,\\
    v_{j}+tv_{j+1} & \text{if }i=j+1.
    \end{cases}
\]
% In matrix form, the action is given by
% \[\sigma_{1} \mapsto\begin{pNiceArray}{cc|c}
%  -t & 1 & \\
% 0 & 1 & \\
% \hline
%  & & I_{n-3}\\
% \end{pNiceArray}, \quad \sigma_{n-1} \mapsto \begin{pNiceArray}{c|cc}
%  I_{n-3} & & \\
%  \hline
%  & 1 & 0 \\
%  & t & -t \\
% \end{pNiceArray}, \quad
% \text{and} \quad
% \sigma_{i} \mapsto\begin{pNiceArray}{c|ccc|c}
%  I_{i-2} & & & & \\
%  \hline
%  & 1 & 0 & 0 & \\
%  & t & -t & 1 & \\
%  & 0 & 0 & 1 \\
%  \hline
%  & & & & I_{n-i-2} \\
% \end{pNiceArray}.
% \]
% By a straightforward calculation of the transposes of the above matrices, we obtain the following:
% \begin{lem}
%     The action of the braid group $\B_{n}$ on $\Bur_{t}^{\inv}$ is given by
% \[
% \sigma_{i}^{-1}(v_{j}) = \begin{cases}
%     v_{j} & i \neq j-1, j, j+1\\
%     t^{-1}v_{j-1}+v_{j} & i = j-1\\
%     -t^{-1}v_{j} & i=j\\
%     v_{j}+v_{j+1} & i=j+1.
%     \end{cases}
% \]
% \end{lem}
% \begin{lem}
%    The action of the braid group $\B_{n}$ on $\Bur_{t}^{\tr}$ is given by
% \[\sigma_{i}^{\tr}(v_{j}) = \begin{cases}
%    v_{j} & i \neq j\\
%    -tv_{1} +v_{2} & i=j=1\\
%    tv_{n-2} - tv_{n-1} & i=j=n-1\\
%    tv_{j-1} - tv_{j} +v_{j+1} & \text{else}.
%    \end{cases}\]
% \end{lem}
% It is known that $\Bur_{t}^{\inv} \cong \Bur_{t^{-1}}$ via the conjugation by the anti-diagonal matrix whose nonzero entries are all $1$. A fortiori, $\Bur_{t}^{\vee} \cong \Bur_{t^{-1}}^{\tr}$.    
For $\zeta\in\bC^{\times}$, recall that $\Bur_{\zeta}$ denotes the specialisation of $\Bur_{t}$ at $t=\zeta$, as defined in~\eqref{eq:def-specialisation}. We first compute the coinvariants of $\Bur_{t}$ and its specialisations.

\begin{prop}\label{prop:coinvariants-Bur}
For $n\ge 3$, $(\Bur_{t})_{\B_{n}} = 0$ and $(\Bur_{\zeta})_{\B_{n}}=0$ for any $\zeta \in \bC^{\times}$. For $n=2$, we have $(\Bur_{t})_{\B_{2}}\cong\bZ[t^{\pm1}]/(1+t)$, and for  $\zeta\in\bC^{\times}$,
\[
(\Bur_{\zeta})_{\B_{2}}\cong
\begin{cases}
0 & \text{if }\zeta\neq-1,\\
\bC & \text{if }\zeta=-1.
\end{cases}
\]
\end{prop}

\begin{proof}
For $n\geq3$, by the definition of the reduced Burau representation, we have
\[
(1-\Bur_{t}(\sigma_{i}))(v_{j})=
\begin{cases}
0 & \text{if }i\notin\{j-1,j,j+1\},\\
-v_{j-1} & \text{if }i=j-1,\\
(1+t)v_{j} & \text{if }i=j,\\
-tv_{j+1} & \text{if }i=j+1.
\end{cases}
\]
For each $1\leq j\leq n-1$, the basis vector $v_{j}$ belongs to $\Image(\Psi_{\Bur_{t}})$ by the second or fourth case above. Therefore, we have
\[
\Image(\Psi_{\Bur_{t}})
=\bZ[t^{\pm1}]\{v_{1},\ldots,v_{n-1}\}
=\Bur_{t}
\]
and so $(\Bur_{t})_{\B_{n}}=0$. Since this argument remains valid under every specialisation $t=\zeta\in\bC^{\times}$ as taking coinvariants commutes with extension of scalars, we also obtain $(\Bur_{\zeta})_{\B_{n}}=0$.

When $n=2$, we have $\Bur_{t}\cong\bZ[t^{\pm1}]$, with the generator of $\B_{2}$ acting by multiplication by $-t$. Hence, $\Image(\Psi_{\Bur_{t}})=(1+t)\bZ[t^{\pm1}]\{v_{1}\}$, which yields the asserted computation of $(\Bur_{t})_{\B_{2}}$. The result for $
(\Bur_{\zeta})_{\B_{2}}$ follows by specialising $t=\zeta\in\bC^{\times}$.
\end{proof}

Applying Lemma~\ref{lem:invariants-coinvariants-duality} to Proposition~\ref{prop:coinvariants-Bur}, we obtain the following corollary.

\begin{coro}\label{coro:invariants-Bur-dual}
For $n\geq 3$, we have $\bigl(\Bur_{t}^{\vee}\bigr)^{\B_{n}}=0$, and $\bigl(\Bur_{\zeta}^{\vee}\bigr)^{\B_{n}}=0$ for each $\zeta\in\bC^{\times}$. For $n=2$, we have $\bigl(\Bur_{t}^{\vee}\bigr)^{\B_{2}}=0$, while, for each $\zeta\in\bC^{\times}$,
\[
\bigl(\Bur_{\zeta}^{\vee}\bigr)^{\B_{2}}\cong
\begin{cases}
0 & \text{if }\zeta\neq-1,\\
\bC & \text{if }\zeta=-1.
\end{cases}
\]
\end{coro}

We now turn to the invariants of $\Bur_{t}$ and its complex specialisations. We begin by computing the coinvariants of the dual representation $\Bur_{t}^{\vee}$. By the relations~\eqref{eq:transpose-dual-invariants-coinvariants}, it suffices to compute the invariants of the transpose representation $\Bur_{t}^{\tr}$. Taking the transpose of the matrices defining the action of the standard generators of $\B_{n}$ on $\Bur_{t}$, we obtain the following formulas: for $n\geq 3$, the action of $\sigma_{i}$ on $\Bur_{t}^{\tr}$ is given by
\[
\Bur_{t}^{\tr}(\sigma_{i})(v_{j})=
\begin{cases}
v_{j} & \text{if }i\neq j,\\
-tv_{1}+v_{2} & \text{if }i=j=1,\\
tv_{n-2}-tv_{n-1} & \text{if }i=j=n-1,\\
tv_{j-1}-tv_{j}+v_{j+1} & \text{if }i=j,\ 2\leq j\leq n-2.
\end{cases}
\]
Using the preceding description of the transpose reduced Burau representation, we obtain the following computation of the coinvariants of its dual.

\begin{prop}\label{prop:coinvariants-BurDual}
For $n\geq 2$, we have
\[
\bigl(\Bur_{t}^{\vee}\bigr)_{\B_{n}}
\cong \bZ[t^{\pm1}]/(1+t+\cdots+t^{n-1}).
\]
Moreover, for $\zeta\in\bC^{\times}$,
\[
\bigl(\Bur_{\zeta}^{\vee}\bigr)_{\B_{n}}
\cong
\begin{cases}
\bC & \text{if }\zeta\neq 1\text{ and }\zeta^{n}=1,\\
0 & \text{otherwise}.
\end{cases}
\]
\end{prop}

\begin{proof}
Assume first that $n\geq 3$. By the preceding formulas, we have
\[
(1-\Bur_{t}^{\tr}(\sigma_{i}))(v_{j})=
\begin{cases}
0 & \text{if }i\neq j,\\
(1+t)v_{1}-v_{2} & \text{if }i=j=1,\\
-tv_{n-2}+(1+t)v_{n-1} & \text{if }i=j=n-1,\\
-tv_{j-1}+(1+t)v_{j}-v_{j+1} & \text{if }i=j,\ 2\leq j\leq n-2.
\end{cases}
\]
We regard the generators of $\Image(\Psi_{\Bur_{t}^{\tr}})$ as the defining relations 
\[
\begin{cases}
(1+t)v_{1}-v_{2}=0,\\
-tv_{j-1}+(1+t)v_{j}-v_{j+1}=0
& \text{for } 2\leq j\leq n-2,\\
-tv_{n-2}+(1+t)v_{n-1}=0
\end{cases}
\]
of the quotient
\[
\bigl(\Bur_{t}^{\tr}\bigr)_{\B_{n}}
=
\frac{\bZ[t^{\pm1}]\{v_{1},\ldots,v_{n-1}\}}
{\Image(\Psi_{\Bur_{t}^{\tr}})}.
\]
We show by induction that $v_{j}=(1+t+\cdots+t^{j-1})v_{1}$ for each $1\leq j\leq n-1$ in the quotient. The base cases follow from the identity $v_{1}=v_{1}$ and the first relation, which gives $v_{2}=(1+t)v_{1}$. Suppose that the formula holds for all $1\leq k\leq j$, for some $2\leq j\leq n-2$. By the second relation, we have
\[
\begin{aligned}
v_{j+1}
&=-tv_{j-1}+(1+t)v_{j}\\
&=-t(1+t+\cdots+t^{j-2})v_{1}
+(1+t)(1+t+\cdots+t^{j-1})v_{1}\\
&=(1+t+\cdots+t^{j})v_{1}.
\end{aligned}
\]
This completes the induction. Hence, each basis element $v_{j}$, for $1\leq j\leq n-1$, is generated by $v_{1}$ in the quotient. Substituting these expressions into the last relation yields $(1+t+\cdots+t^{n-1})v_{1}=0$. Therefore, the quotient is generated by $v_{1}$, subject to this single relation. More precisely, we have
\[
\begin{aligned}
\bigl(\Bur_{t}^{\tr}\bigr)_{\B_{n}}
&\cong
\frac{\bZ[t^{\pm1}]\{v_{1},\ldots,v_{n-1}\}}
{\bigl((1+t+\cdots+t^{n-1})v_{1},\;
v_j-(1+t+\cdots+t^{j-1})v_{1}:2\leq j\leq n-1\bigr)}
\\
&\cong
\frac{\bZ[t^{\pm1}]\{v_{1}\}}
{\bigl((1+t+\cdots+t^{n-1})v_{1}\bigr)}
\cong \bZ[t^{\pm1}]/(1+t+\cdots+t^{n-1}).
\end{aligned}
\]

We now consider the complex specialisations $\Bur_{\zeta}^{\tr}$, for $\zeta\in\bC^{\times}$. Since taking coinvariants commutes with extension of scalars, we deduce that
\[
\bigl(\Bur_{\zeta}^{\tr}\bigr)_{\B_{n}}
\cong \bC/(1+\zeta+\cdots+\zeta^{n-1}).
\]
If $\zeta=1$, then $1+\zeta+\cdots+\zeta^{n-1}=n\neq 0$. If $\zeta\neq 1$, we have
\[1+\zeta+\cdots+\zeta^{n-1} =\frac{1-\zeta^{n}}{1-\zeta}\]
which vanishes precisely when $\zeta^{n}=1$. Since every nonzero complex number is invertible, this yields the desired result for $\bigl(\Bur_{\zeta}^{\tr}\bigr)_{\B_{n}}$.

Finally, consider the case $n=2$. Since $\Bur_{t}$ is one-dimensional, with the generator $\sigma_{1}$ of $\B_{2}$ acting by multiplication by $-t$, we have $\Bur_{t}^{\tr}=\Bur_{t}$. Hence, we conclude that $\bigl(\Bur_{t}^{\tr}\bigr)_{\B_{2}}
\cong \bZ[t^{\pm1}]/(1+t)$. Similarly, for each $\zeta\in\bC^{\times}$, we have
\[
\bigl(\Bur_{\zeta}^{\tr}\bigr)_{\B_{2}}
\cong
\begin{cases}
\bC & \text{if }\zeta=-1,\\
0 & \text{otherwise}.
\end{cases}
\]
Thus, the asserted formulas also hold for $n=2$.
\end{proof}

Combining Proposition~\ref{prop:coinvariants-BurDual} with Lemma~\ref{lem:invariants-coinvariants-duality}, we obtain the following computation of the invariants of the reduced Burau representation and its complex specialisations.

\begin{coro}\label{coro:invariants-Bur}
For $n\ge 2$, $\Bur_{t}^{\B_{n}} = 0$, and for $\zeta\in\bC^{\times}$,
\[
\Bur_{\zeta}^{\B_{n}}\cong
\begin{cases}
\bC & \text{if }\zeta\neq 1\text{ and }\zeta^n=1,\\
0 & \text{otherwise}.
\end{cases}
\]
\end{coro}

%-------------------------------------------------------

\section{The Lawrence--Krammer--Bigelow representation}\label{s:LKB-representation}

Having studied the low-dimensional representations of $\B_{n}$, we now turn to the Lawrence--Krammer--Bigelow representation. In this section, we compute its invariants and coinvariants, together with those of its dual.

%-------------------------------------------------------

\subsection{The Lawrence--Krammer--Bigelow representation}\label{ss:LKB-representation}

We begin by recalling the explicit formulas defining the Lawrence--Krammer--Bigelow representation, which will serve as the starting point for our computations. The Lawrence--Krammer--Bigelow representation $\LKB_{t,q}$ is a representation of the braid group $\B_{n}$ on a free $\bZ[t^{\pm1},q^{\pm1}]$-module of rank $n(n-1)/2$, arising from the homological representations introduced by Lawrence. Its explicit action was computed by Bigelow~\cite[Th.~4.1]{bigelow2001braid}, with a sign correction recorded in~\cite[\S 5.2]{bigelowthelkrepresentation}. An alternative description, using a different basis and parameter conventions, is given by Krammer~\cite[\S 3]{KrammerLK}.

\begin{prop}[Bigelow \cite{bigelow2001braid,bigelowthelkrepresentation}]\label{prop:formula-lkb}
The representation $\LKB_{t,q}$ has a basis $\{v_{j,k}:1\leq j<k\leq n\}$, with the action of the standard generators $\sigma_i$ of $\B_{n}$ given by
\begin{equation}\label{eq:lkb}
    \LKB_{t,q}(\sigma_{i})(v_{j,k}) = \begin{cases}
    v_{j,k} & i \neq j-1, j, k-1, k \quad (a)\\
    qv_{j-1,k} + (q^2-q)v_{j-1,j} + (1-q)v_{j,k} & i=j-1 \hfill (b)\\
    v_{j+1,k} & i=j<k-1 \hfill (c)\\
    qv_{j,k-1} + (1-q)v_{j,k}-t(q^2-q)v_{k-1,k} & i=k-1>j \hfill (d)\\
    v_{j,k+1} & i=k \hfill (e)\\
    -tq^2 v_{j,k} = -tq^2 v_{j,j+1} & i=j=k-1. \hfill (f)
\end{cases}
\end{equation}
\end{prop}

We determine the explicit action of the transpose of the Lawrence--Krammer--Bigelow representation.

\begin{lem}\label{lem:lkb-transpose}
The action of the standard generators $\sigma_i$ of $\B_{n}$ on $\LKB_{t,q}^{\tr}$ is given by
\begin{equation}\label{eq:lkb-transpose}
\LKB_{t,q}^{\tr}(\sigma_{i})(v_{j,k}) = \begin{cases}
    v_{j,k} & i \neq j-1, j, k-1, k \quad (a)\\
    v_{j-1,k}+(1-q)v_{j,k} & i=j-1 \hfill (b)\\
    qv_{j+1,k} & i=j<k-1 \hfill (c)\\
    v_{j,k-1}+(1-q)v_{j,k} & i=k-1>j \hfill (d)\\
    qv_{j,k+1} & i=k \hfill (e)\\
    \displaystyle -tq^2v_{j,j+1}+(q^2-q)\sum_{h=j+2}^{n} v_{j+1,h} + t(q-q^2)\sum_{h=1}^{j-1} v_{h,j+1} & i=j=k-1. \hfill (f)
    \end{cases}
\end{equation}
\end{lem}

\begin{proof}
    To compute the formula for the action of $\B_{n}$ on the transpose $\LKB_{t,q}^{\tr}$, it is best to first fix $1 \le i \le n-1$ and rewrite Equation~\ref{eq:lkb} in terms of $i$:
    \begin{equation}\label{eq:lkb-i}
    \LKB_{t,q}(\sigma_{i})(v_{j,k}) = \begin{cases}
    v_{j,k} & j,k \ne i,i+1 \hspace{1.5cm} (a)\\
    qv_{i,k} + (q^2-q)v_{i,i+1} + (1-q)v_{i+1,k} & j=i+1; k\ne i,i+1 \hfill (b)\\
    v_{i+1,k} & j=i;k\ne i,i+1 \hfill (c)\\
    qv_{j,i} + (1-q)v_{j,i+1}-t(q^2-q) v_{i,i+1} & k=i+1; j\ne i,i+1 \hfill (d)\\
    v_{j,i+1} & i=k; j \ne i,i+1 \hfill (e)\\
    -tq^2 v_{i,i+1} & j=i; k=i+1. \hfill (f)
    \end{cases}
    \end{equation}
    For each $1 \le j < k \le n$, we want to find all basis elements $v$ such that the coefficient of $v_{j,k}$ in $\LKB_{t,q}(\sigma_{i})(v)$ is not zero. This number is precisely the coefficient of $v$ in $\LKB_{t,q}^{\tr}(\sigma_{i})(v_{j,k})$. For example, if $j = i+1$ and $k \ne i, i+1$, the term $v_{i+1,k}$ appears in $\LKB_{t,q}(\sigma_{i})(v_{i+1,k})$ with coefficient $1-q$ and in $\LKB_{t,q}(\sigma_{i})(v_{i,k})$ with coefficient $1$. Thus $\LKB_{t,q}^{\tr}(\sigma_{i})(v_{i+1,k}) = (1-q)v_{i+1,k}+v_{i,k}$. Repeating this calculation for all cases, we have
    \begin{equation}\label{eq:lkb-transpose_{i}}
    \LKB_{t,q}^{\tr}(\sigma_{i})(v_{j,k}) = \begin{cases}
        v_{j,k} & j,k \ne i,i+1 \hspace{1.5cm} (a)\\
        v_{i,k}+(1-q)v_{i+1,k} & j=i+1; k\ne i,i+1 \hfill (b)\\
        qv_{i+1,k} & j=i;k\ne i,i+1 \hfill (c)\\
        v_{j,i}+(1-q)v_{j,i+1} & k=i+1; j\ne i,i+1 \hfill (d)\\
        qv_{j,i+1} & i=k; j \ne i,i+1 \hfill (e)\\
        \displaystyle -tq^2v_{i,i+1}+(q^2-q)\sum_{h=i+2}^{n} v_{i+1,h} + t(q-q^2)\sum_{h=1}^{i-1} v_{h,i+1} & j=i; k=i+1. \hfill (f)
    \end{cases}
    \end{equation}
    The last formula represents the most interesting case. Here, the term $v_{i,i+1}$ appears in $\LKB_{t,q}(\sigma_{i})(v_{i,i+1})$ with coefficient $-tq^2$, in $\LKB_{t,q}(\sigma_{i})(v_{i+1,k})$ with coefficient $q^2-q$ for all $k \ge i+2$, and in $\LKB_{t,q}^{\tr}(\sigma_{i})(v_{j,i+1})$ with coefficient $-t(q^2-q)$ for all $j \le i-1$. By reindexing Equation~\ref{eq:lkb-transpose_{i}} in terms of $j$ and $k$, we obtain Equation~\ref{eq:lkb-transpose}.
\end{proof}

%--------------------------------------------------

\subsection{Coinvariants}

We begin by computing the coinvariants of the Lawrence--Krammer--Bigelow representation $\LKB_{t,q}$ over its natural coefficient ring $\bZ[t^{\pm1},q^{\pm1}]$. Using the explicit action of the standard braid generators given in Proposition~\ref{prop:formula-lkb}, we determine the relations defining the coinvariant module and obtain the following result.

For $\bzeta=(\zeta_{t},\zeta_{q})\in(\bC^{\times})^{2}$, we denote by $\LKB_{\bzeta}$ the complex specialisation of $\LKB_{t,q}$ obtained by setting $t=\zeta_{t}$ and $q=\zeta_{q}$, as defined in~\eqref{eq:def-specialisation}.

\begin{prop}\label{prop:coinvariants-LKB}
For $n\geq 3$, we have
\[
(\LKB_{t,q})_{\B_{n}}
\cong
\bZ[t^{\pm1},q^{\pm1}]/(1+t,1-q).
\]
Moreover, for $\bzeta=(\zeta_{t},\zeta_{q})\in(\bC^{\times})^{2}$ and $n\geq 3$, we have
\[
(\LKB_{\bzeta})_{\B_{n}}
\cong
\begin{cases}
\bC & \text{if }\zeta_{t}=-1\text{ and }\zeta_{q}=1,\\
0 & \text{otherwise}.
\end{cases}
\]
For $n=2$, we have instead
\[
(\LKB_{t,q})_{\B_{2}}
\cong
\bZ[t^{\pm1},q^{\pm1}]/(1+tq^{2}),
\]
and
\[
(\LKB_{\bzeta})_{\B_{2}}
\cong
\begin{cases}
\bC & \text{if }\zeta_{t}\zeta_{q}^{2}=-1,\\
0 & \text{otherwise}.
\end{cases}
\]
\end{prop}

\begin{proof}
Assume first that $n\geq 3$. By the explicit formulas of Proposition~\ref{prop:formula-lkb}, for $1\leq i\leq n-1$ and $1\leq j<k\leq n$, we have
\begin{equation}
\begin{split}
&(1 - \LKB_{t,q}(\sigma_{i}))(v_{j,k})\\
& = \begin{cases}
v_{j,k} - v_{j,k} = 0 & i \neq j-1, j,k-1,k \quad (a)\\
v_{j,k} - (qv_{j-1,k} + (q^2-q)v_{j-1,j} + (1-q)v_{j,k}) & i=j-1 \hfill (b)\\
v_{j,k} - v_{j+1,k} & i=j<k-1 \hfill (c)\\
v_{j,k} - (qv_{j,k-1} + (1-q)v_{j,k}-t(q^2-q)v_{k-1,k}) & i=k-1>j \hfill (d)\\
v_{j,k} - v_{j,k+1} & i=k \hfill (e)\\
v_{j,j+1} - (-tq^2 v_{j,j+1}) & i=j=k-1. \hfill (f)
\end{cases}
\end{split}
\end{equation}
We regard these generators of $\Image(\Psi_{\LKB_{t,q}})$ as the defining relations of the quotient
\[
(\LKB_{t,q})_{\B_{n}} = \frac{\bZ[t^{\pm1},q^{\pm1}] \{v_{j,k}:1\leq j<k\leq n\}} {\Image(\Psi_{\LKB_{t,q}})}.
\]
More precisely, the defining relations are given by the following system:
\begin{equation}\label{eq:coinvariants-LKB-system}
\begin{cases}
qv_{j,k} - qv_{j-1,k} + (q-q^2)v_{j-1,j} = 0 & 2 \le j=i+1 < k \le n \hspace{1cm} (b)\\
v_{j,k} - v_{j+1,k} =0 & 
1 \le j=i<k-1 \le n-1 \hfill (c)\\
qv_{j,k} - qv_{j,k-1} +t(q^2-q)v_{k-1,k} =0 & 1 \le j < i=k-1 < k \le n \hfill (d)\\
v_{j,k} - v_{j,k+1} =0 & 1 \le j < k=i \le n-1 \hfill (e)\\
(1+tq^2 )v_{j,j+1} =0 & 1\le j=i=k-1 \le n-1. \hfill (f)
\end{cases}
\end{equation}
By~\eqref{eq:coinvariants-LKB-system}(e), we have $v_{1,k}=v_{1,2}$ for each $2\leq k\leq n$ in the quotient. Similarly, \eqref{eq:coinvariants-LKB-system}(c) implies that $v_{j,k}=v_{1,k}$ for each $1\leq j<k\leq n$. Therefore, all the basis elements are identified with $v_{1,2}$ in the quotient. Under these identifications, the relations \eqref{eq:coinvariants-LKB-system}(b) and \eqref{eq:coinvariants-LKB-system}(d) both reduce to $(1-q)v_{1,2}=0$, since $q$ and $t$ are invertible. Moreover, \eqref{eq:coinvariants-LKB-system}(f) reduces to $(1+tq^{2})v_{1,2}=0$. Since $(1+t)v_{1,2} = (1+tq^{2})v_{1,2} +t(1-q^{2})v_{1,2} =0$, the latter relation, together with $(1-q)v_{1,2}=0$, is equivalent to $(1+t)v_{1,2}=0$. Therefore, the quotient is generated by $v_{1,2}$, subject precisely to the relations $(1-q)v_{1,2}=0$ and $(1+t)v_{1,2}=0$.  We conclude that
\[
\begin{aligned}
(\LKB_{t,q})_{\B_{n}}
&\cong
\frac{\bZ[t^{\pm1},q^{\pm1}]\{v_{1,2}\}}
{((1-q)v_{1,2},(1+t)v_{1,2})}
&\cong
\bZ[t^{\pm1},q^{\pm1}]/(1-q,1+t).
\end{aligned}
\]
When $n=2$, the representation $\LKB_{t,q}$ is free of rank one, generated by $v_{1,2}$, and $\Psi_{\LKB_{t,q}}(v_{1,2})=(1+tq^{2})v_{1,2}$. We deduce that $(\LKB_{t,q})_{\B_{2}} \cong \bZ[t^{\pm1},q^{\pm1}]/(1+tq^{2})$.

We now consider the complex specialisations of $\LKB_{t,q}$. For $n\geq 3$, since taking coinvariants commutes with extension of scalars, we deduce that
\[
(\LKB_{\bzeta})_{\B_{n}}
\cong
\bC/(1+\zeta_{t},1-\zeta_{q}).
\]
Since $\bC$ is a field, this quotient is isomorphic to $\bC$ if both $\zeta_{t}=-1$ and $\zeta_{q}=1$, and vanishes otherwise. Similarly, when $n=2$, we obtain
\[
(\LKB_{\bzeta})_{\B_{2}}
\cong
\bC/(1+\zeta_{t}\zeta_{q}^{2}).
\]
Thus, the coinvariants are one-dimensional if $\zeta_{t}\zeta_{q}^{2}=-1$, and vanish otherwise.
\end{proof}

Since the coinvariants computed in Proposition~\ref{prop:coinvariants-LKB} are torsion modules over the integral domain $\bZ[t^{\pm1},q^{\pm1}]$, Lemma~\ref{lem:invariants-coinvariants-duality} immediately yields the following description of the invariants of the dual Lawrence--Krammer--Bigelow representation and its complex specialisations.

\begin{coro}\label{coro:invariants-LKB-dual}
For $n\geq 2$, we have $(\LKB_{t,q}^{\vee})^{\B_{n}}=0$. Moreover, for $\bzeta=(\zeta_{t},\zeta_{q})\in(\bC^{\times})^{2}$ and $n\geq 3$, we have
\[
(\LKB_{\bzeta}^{\vee})^{\B_{n}}
\cong
\begin{cases}
\bC & \text{if }\zeta_{t}=-1\text{ and }\zeta_{q}=1,\\
0 & \text{otherwise}.
\end{cases}
\]
For $n=2$, we have
\[
(\LKB_{\bzeta}^{\vee})^{\B_{2}}
\cong
\begin{cases}
\bC & \text{if }\zeta_{t}\zeta_{q}^{2}=-1,\\
0 & \text{otherwise}.
\end{cases}
\]
\end{coro}

%--------------------------------------------------

\subsection{Invariants}

We now turn to the invariants of the Lawrence--Krammer--Bigelow representation $\LKB_{t,q}$. As in the preceding sections, our approach relies on the duality between invariants and coinvariants established in Lemma~\ref{lem:invariants-coinvariants-duality}. We therefore begin by computing the coinvariants of the dual representation $\LKB_{t,q}^{\vee}$, from which we will deduce the invariants of $\LKB_{t,q}$.

\begin{prop}\label{prop:coinvariants-LKB-dual}
For $n\geq 2$, we have
\[
(\LKB_{t,q}^{\vee})_{\B_{n}}
\cong
\bZ[t^{\pm1},q^{\pm1}]/(1+tq^{n}).
\]
Moreover, for $\bzeta=(\zeta_{t},\zeta_{q})\in(\bC^{\times})^{2}$, we have
\[
(\LKB_{\bzeta}^{\vee})_{\B_{n}}
\cong
\begin{cases}
\bC & \text{if }\zeta_{t}\zeta_{q}^{n}=-1,\\
0 & \text{otherwise}.
\end{cases}
\]
\end{prop}

\begin{proof}
By the relations~\eqref{eq:transpose-dual-invariants-coinvariants}, it suffices to compute the coinvariants of the transpose representation $\LKB_{t,q}^{\tr}$. Assume first that $n\geq 3$. For $1\leq i\leq n-1$ and $1\leq j<k\leq n$, we have
\begin{equation}
\begin{split}
&(1-\LKB_{t,q}^{\tr}(\sigma_{i}))(v_{j,k})\\
&=
\begin{cases}
0
& i\notin\{j-1,j,k-1,k\} \quad (a),\\
v_{j,k}-(v_{j-1,k}+(1-q)v_{j,k})
& i=j-1 \hfill (b),\\
v_{j,k}-qv_{j+1,k}
& i=j<k-1 \hfill (c),\\
v_{j,k}-(v_{j,k-1}+(1-q)v_{j,k})
& i=k-1>j \hfill (d),\\
v_{j,k}-qv_{j,k+1}
& i=k \hfill (e),\\
v_{j,j+1}-\left(
\displaystyle -tq^{2}v_{j,j+1}
+(q^{2}-q)\sum_{h=j+2}^{n}v_{j+1,h}
+t(q-q^{2})\sum_{h=1}^{j-1}v_{h,j+1}
\right)
& i=j=k-1. \hfill (f)
\end{cases}
\end{split}
\end{equation}
We regard these generators of $\Image(\Psi_{\LKB_{t,q}^{\tr}})$ as the defining relations of the quotient
\[
(\LKB_{t,q}^{\tr})_{\B_{n}}
=
\frac{\bZ[t^{\pm1},q^{\pm1}]\{v_{j,k}:1\leq j<k\leq n\}}{\Image(\Psi_{\LKB_{t,q}^{\tr}})}.
\]
Since relations (c) and (e) are reindexed versions of
relations (b) and (d), respectively, the defining
relations reduce to the following system:
\begin{equation}\label{eq:coinvariants-LKBdual-system}
\begin{cases}
qv_{j,k} - v_{j-1,k} =0 & 2 \le j < k \le n \hspace{1.6cm} (b)\\
qv_{j,k} - v_{j,k-1} =0 & 1 \le j < k-1 \le n-1 \hfill (d)\\
\displaystyle (1+tq^2)v_{j,j+1}+(q-q^2)\sum_{h=j+2}^{n} v_{j+1,h} + t(q^2-q)\sum_{h=1}^{j-1} v_{h,j+1} =0 & 1\le j=k-1 \le n-1. \hfill (f)
\end{cases}
\end{equation}
By~\eqref{eq:coinvariants-LKBdual-system}(d), we have $v_{1,k}=q^{-k+2}v_{1,2}$ for every $2\leq k\leq n$ in the quotient. Similarly, \eqref{eq:coinvariants-LKBdual-system}(b) implies that $v_{j,k}=q^{-j+1}v_{1,k}$ for every $1\leq j<k\leq n$, and so $v_{j,k}=q^{3-j-k}v_{1,2}$. Thus, the quotient is generated by $v_{1,2}$. Under these identifications, the final relation \eqref{eq:coinvariants-LKBdual-system}(f) simplifies to
\[\displaystyle \left( q^{2-2j}+tq^{4-2j} + q^{3-j-n}(1-q)\sum_{h=0}^{n-j-2} q^{h} + tq^{4-2j}(q-1)\sum_{h=0}^{j-2} q^{h} \right) v_{1,2} =0.
\]
Using the geometric series identities, this simplifies to
\[\displaystyle \left( q^{2-2j}+tq^{4-2j} + q^{3-j-n}(1-q^{n-j-1})+ tq^{4-2j}(q^{j-1}-1) \right) v_{1,2} = q^{3-j-n}(1+tq^{n})v_{1,2} =0.
\]
Since $q$ is invertible, all these relations are equivalent to
\[
(1+tq^{n})v_{1,2}=0.
\]
Therefore, the quotient is generated by $v_{1,2}$, subject precisely to this single relation. We conclude that
\[
(\LKB^{\vee}_{t,q})_{\B_{n}} \cong (\LKB^{\tr}_{t,q})_{\B_{n}} \cong \dfrac{\bZ[t^{\pm1},q^{\pm1}]\{v_{j,k}:1 \le j < k \le n\}}{\Image(\Psi_{\LKB_{t,q}^{\tr}})} \cong \dfrac{\bZ[t^{\pm1},q^{\pm1}]\{v_{1,2}\}}{((1+tq^{n})v_{1,2})} \cong \dfrac{\bZ[t^{\pm1},q^{\pm1}]}{(1+tq^{n})}.
\]

When $n=2$, the representation $\LKB_{t,q}^{\tr}$ is free of rank one, generated by $v_{1,2}$, and $\Psi_{\LKB_{t,q}^{\tr}}(v_{1,2}) =(1+tq^{2})v_{1,2}$. Hence, we have $(\LKB_{t,q}^{\vee})_{\B_{2}} \cong \bZ[t^{\pm1},q^{\pm1}]/(1+tq^{2})$, which extends the preceding formula to $n=2$.

We now consider the complex specialisations of the dual Lawrence--Krammer--Bigelow representation. Since taking coinvariants commutes with extension of scalars, we deduce that $(\LKB_{\bzeta}^{\vee})_{\B_{n}} \cong \bC/(1+\zeta_{t}\zeta_{q}^{n})$. Since $\bC$ is a field, this quotient is one-dimensional if $\zeta_{t}\zeta_{q}^{n}=-1$, and vanishes otherwise.
\end{proof}

Since the coinvariants computed in Proposition~\ref{prop:coinvariants-LKB-dual} are torsion modules over the integral domain $\bZ[t^{\pm1},q^{\pm1}]$, Lemma~\ref{lem:invariants-coinvariants-duality} immediately yields the following computation of the invariants of the Lawrence--Krammer--Bigelow representation and its complex specialisations.

\begin{coro}\label{coro:invariants-LKB}
For $n\geq 2$, we have
\[
(\LKB_{t,q})^{\B_{n}}=0.
\]
Moreover, for $\bzeta=(\zeta_{t},\zeta_{q})\in(\bC^{\times})^{2}$, we have
\[
(\LKB_{\bzeta})^{\B_{n}}
\cong
\begin{cases}
\bC & \text{if }\zeta_{t}\zeta_{q}^{n}=-1,\\
0 & \text{otherwise}.
\end{cases}
\]
\end{coro}

%--------------------------------------------------

\section{Tensor powers of the Tong--Yang--Ma representation} \label{s:tensor-TYM}

In this section, we study tensor products of Tong--Yang--Ma representations and compute their coinvariants and invariants. Throughout this section, we write $R_{d}:=\bZ[t_{1}^{\pm1},\dots,t_{d}^{\pm1}]$.

%--------------------------------------------------

\subsection{Tensor power decomposition}

Let $d\geq 2$. We consider the tensor product
\[
V_{d}:=\TYM_{t_{1}}\otimes\cdots\otimes\TYM_{t_{d}}
\]
over $R_{d}$, where each factor is extended to $R_{d}$. As before, let $v_{j}$ denote the $j^{\mathrm{th}}$ standard basis vector of $R_{d}^{\oplus n}$. We begin by observing that the action of the braid generators on $V_{d}$ has a particularly simple form. Let $(i\,\,i+1)$ denote the transposition of $\{1,\dots,n\}$ exchanging $i$ and $i+1$. For any $1\leq j_{1},\dots,j_{d}\leq n$, we have
\begin{equation}\label{eq:TYM-tensor-action}
V_{d}(\sigma_{i})(v_{j_{1}}\otimes\cdots\otimes v_{j_{d}})
=
\left(\prod_{\ell:\,j_{\ell}=i}t_{\ell}\right)
v_{(i\,\,i+1)(j_{1})}\otimes\cdots\otimes
v_{(i\,\,i+1)(j_{d})}.
\end{equation}
In particular, each braid generator acts by permuting the standard basis elements of $V_{d}$, up to multiplication by units in $R_{d}$. We thus distinguish two natural subrepresentations of $V_{d}$.

\begin{lem}\label{lem:TYM-tensor-iag}
The $R_{d}$-submodule of $V_{d}$ generated by $\{v_{j}\otimes\cdots\otimes v_{j} : 1\leq j\leq n\}$ is a $\B_{n}$-subrepresentation, called the \emph{diagonal subrepresentation} $\Delta(V_{d})$. Moreover, there is an isomorphism of $\B_{n}$-representations
\[
\Delta(V_{d})\cong\TYM_{t_{1}\cdots t_{d}}.
\]
\end{lem}

\begin{proof}
By~\eqref{eq:TYM-tensor-action}, the action of each standard generator $\sigma_{i}$ on the generators of $\Delta(V_{d})$ is given by
\[
V_{d}(\sigma_{i})
(v_{j}\otimes\cdots\otimes v_{j})
=
\begin{cases}
v_{j}\otimes\cdots\otimes v_{j}
& i\notin\{j-1,j\},\\
v_{j-1}\otimes\cdots\otimes v_{j-1}
& i=j-1,\\
(t_{1}\cdots t_{d})
(v_{j+1}\otimes\cdots\otimes v_{j+1})
& i=j.
\end{cases}
\]
Thus, $\Delta(V_{d})$ is preserved by the braid-group action. Furthermore, these formulas agree with those defining the Tong--Yang--Ma representation $\TYM_{t_{1}\cdots t_{d}}$. Hence, the $R_{d}$-linear map $\Delta(V_{d})\to\TYM_{t_{1}\cdots t_{d}}$ sending $v_{j}\otimes\cdots\otimes v_{j}$ to $v_{j}$ for each $1\leq j\leq n$ is a $\B_{n}$-equivariant isomorphism.
\end{proof}

\begin{lem}\label{lem:TYMn_aDiag}
The $R_{d}$-submodule of $V_{d}$ generated by $\{v_{j_{1}}\otimes\cdots\otimes v_{j_{d}} : 1\leq j_{1},\dots,j_{d}\leq n,\; j_{k}\neq j_{\ell}\text{ whenever }k\neq\ell\}$ is a $\B_{n}$-subrepresentation, called the \emph{anti-diagonal subrepresentation} $\aDiag(V_{d})$.
\end{lem}
\begin{proof}
By~\eqref{eq:TYM-tensor-action}, each standard generator $\sigma_{i}$ acts on the basis elements $v_{j_{1}}\otimes\cdots\otimes v_{j_{d}}$ by simultaneously applying the transposition $(i\,\,i+1)$ to all their indices, up to multiplication by a unit in $R_{d}$. Since this transposition preserves the property of having pairwise distinct indices, it induces a permutation of the basis of $\aDiag(V_{d})$. Hence, each $\sigma_{i}$ restricts to an automorphism of $\aDiag(V_{d})$, which is therefore a $\B_{n}$-subrepresentation of $V_{d}$.
\end{proof}

The diagonal and anti-diagonal subrepresentations are two extreme cases of a more general decomposition of $V_{d}$. Let $\cP_{d}$ denote the set of partitions of $[d]:=\{1,\dots,d\}$, i.e.~each element $\{p_{1},p_{2},\dots,p_{m}\}$ of $\cP_{d}$ consists of disjoint subsets $p_{i}$ of $[d]$ such that $\cup_{i} p_{i} = [d]$. We write $\cP_{d}^{n}$ for the subset of partitions with at most $n$ parts. In particular, $\cP_{d}^{n}=\cP_{d}$ whenever $n\geq d$. For each nonempty subset $p_{i}\subseteq[d]$, we define
\[
t_{p_{i}}:=\prod_{k\in p_{i}}t_{k}.
\]
Whenever we write a partition as $p=\{p_{1},\dots,p_{m}\}$, we order its parts by their least elements. The following decomposition expresses $V_{d}$ as a direct sum of anti-diagonal subrepresentations associated with the partitions of $[d]$.

\begin{lem}\label{lem:TYM-tensor-decomposition}
There is an isomorphism of $\B_{n}$-representations
\[
V_{d}
\cong
\bigoplus_{\{p_{1},\dots,p_{m}\}\in\cP_{d}^{n}}
\aDiag\left(\TYM_{t_{p_{1}}}\otimes\cdots\otimes
\TYM_{t_{p_{m}}}\right).
\]
\end{lem}

\begin{proof}
Recall that $V_{d}$ has the standard $R_{d}$-basis $\{v_{j_{1}}\otimes\cdots\otimes v_{j_{d}} : 1\leq j_{1},\dots,j_{d}\leq n\}$. Each such basis element determines a unique partition $p=\{p_{1},\dots,p_{m}\}$ of $[d]$, characterised by the condition that $j_{k}=j_{\ell}$ if and only if $k$ and $\ell$ belong to the same part of $p$. We say that the basis element is of \emph{type $p$}. For each $p\in\cP_{d}^{n}$, let $V_{p}$ denote the $R_{d}$-submodule of $V_{d}$ generated by all basis elements of type $p$. By~\eqref{eq:TYM-tensor-action}, the action of each standard braid generator simultaneously permutes the indices $j_{1},\dots,j_{d}$, without changing which indices are equal. Hence, the action preserves the type of each basis element, and $V_{p}$ is a $\B_{n}$-subrepresentation of $V_{d}$. Since every basis element has a unique type, and a partition with more than $n$ parts cannot occur, we obtain a direct sum decomposition
\begin{equation}\label{eq:decomposition-V-d-direct-sum}
V_{d}\cong\bigoplus_{p\in\cP_{d}^{n}}V_{p}.    
\end{equation}
It remains to identify the summands. Fix a partition $p=\{p_{1},\dots,p_{m}\}\in\cP_{d}^{n}$ and consider the $R_{d}$-linear map
\[
\begin{aligned}
V_{p}
&\longrightarrow
\aDiag\left(
\TYM_{t_{p_{1}}}\otimes\cdots\otimes
\TYM_{t_{p_{m}}}
\right),\\
v_{j_{1}}\otimes\cdots\otimes v_{j_{d}}
&\longmapsto
v_{j_{\min(p_{1})}}\otimes\cdots\otimes
v_{j_{\min(p_{m})}}.
\end{aligned}
\]
This map is well-defined because the indices $j_{k}$ are equal within each part of $p$ and distinct between different parts. Moreover, it defines a bijection between the standard bases of the source and target, and hence an isomorphism of $R_{d}$-modules. To verify the $\B_{n}$-equivariance, we observe that whenever the indices corresponding to a part $p_{i}$ are equal to $j$, the action of $\sigma_{j}$ on the corresponding tensor factors contributes the scalar
\[
\prod_{k\in p_{i}}t_{k}=t_{p_{i}}.
\]
Therefore, the action on $V_{p}$ agrees, under the preceding identification, with the action on the anti-diagonal subrepresentation of
$\TYM_{t_{p_{1}}}\otimes\cdots\otimes\TYM_{t_{p_{m}}}$, as described by~\eqref{eq:TYM-tensor-action}. Thus, the map is $\B_{n}$-equivariant, and we obtain
\[
V_{p}
\cong
\aDiag\left(
\TYM_{t_{p_{1}}}\otimes\cdots\otimes
\TYM_{t_{p_{m}}}
\right).
\]
Combining these isomorphisms with the decomposition \eqref{eq:decomposition-V-d-direct-sum} completes the proof.
\end{proof}

\subsection{Coinvariants and invariants}\label{ssec:tensor_tym}

We now study the coinvariants and invariants of tensor products of Tong--Yang--Ma representations. We begin with the coinvariants of the anti-diagonal subrepresentation, which will serve as the key ingredient in our subsequent computations.

\begin{lem}\label{lem:coinvariants-TYMn-aDiag}
For $n\geq d+1$, there is an isomorphism of $R_{d}$-modules
\[
\bigl(\aDiag(V_{d})\bigr)_{\B_{n}}
\cong R_{d}/(1-t_{1},\dots,1-t_{d}).
\]
\end{lem}

\begin{proof}
Recall that $\aDiag(V_{d})$ has the $R_{d}$-basis consisting of the elements $v_{j_{1},\dots,j_{d}}
:=v_{j_{1}}\otimes\cdots\otimes v_{j_{d}}$, where $1\leq j_{1},\dots,j_{d}\leq n$ are pairwise distinct. We regard the elements $(1-V_{d}(\sigma_{i}))(v_{j_{1},\dots,j_{d}})$ as the defining relations of the coinvariant module $\aDiag(V_{d})_{\B_{n}}$. By~\eqref{eq:TYM-tensor-action}, each standard generator $\sigma_{i}$ acts on the basis elements by simultaneously applying the transposition $(i\,\,i+1)$ to their indices, up to multiplication by a unit in $R_{d}$. Since $\Sym_{n}$ acts transitively on the set of ordered $d$-tuples of pairwise distinct elements of $[n]$, every basis element is identified in the coinvariant module with a unit multiple of $v_{1,\dots,d}$. We deduce that $\aDiag(V_{d})_{\B_{n}}$ is generated as an $R_{d}$-module by the class of $v_{1,\dots,d}$.

We now determine the relations satisfied by this generator. Fix $1\leq k\leq d$. Since $n\geq d+1$, we may choose a basis element $w_{k}$ whose $k^{\mathrm{th}}$ index is $d$, while its other indices are precisely $1,\dots,d-1$, in any order. In particular, none of its indices is equal to $d+1$. Applying~\eqref{eq:TYM-tensor-action} twice, we obtain $V_{d}(\sigma_{d}^{2})(w_{k})=t_{k}w_{k}$. Since $\sigma_{d}^{2}$ acts trivially on the coinvariants, it follows that $(1-t_{k})w_{k}=0$ in $\aDiag(V_{d})_{\B_{n}}$. As $w_{k}$ is identified with a unit multiple of $v_{1,\dots,d}$, we deduce that
\[
(1-t_{k})v_{1,\dots,d}=0.
\]
This holds for every $1\leq k\leq d$. Hence, there is a surjective $R_{d}$-module homomorphism
\[
R_{d}/(1-t_{1},\dots,1-t_{d})
\twoheadrightarrow
\aDiag(V_{d})_{\B_{n}}
\]
sending the class of $1$ to the class of $v_{1,\dots,d}$. It remains to show that there are no additional relations. Consider the $R_{d}$-linear map
\[
\varepsilon\colon
\aDiag(V_{d})\longrightarrow
R_{d}/(1-t_{1},\dots,1-t_{d})
\]
sending every standard basis element $v_{j_{1},\dots,j_{d}}$ to the class of $1$. By~\eqref{eq:TYM-tensor-action}, the action of each $\sigma_{i}$ permutes these basis elements up to multiplication by a monomial in the variables $t_{1},\dots,t_{d}$. Since all these monomials become $1$ in the target, the map $\varepsilon$ is $\B_{n}$-invariant. It thus induces a surjective $R_{d}$-module homomorphism
\[
\aDiag(V_{d})_{\B_{n}}
\twoheadrightarrow
R_{d}/(1-t_{1},\dots,1-t_{d}).
\]
The composite of the two surjective homomorphisms constructed above is the identity on $R_{d}/(1-t_{1},\dots,1-t_{d})$. They are therefore inverse isomorphisms, proving that
\[
\aDiag(V_{d})_{\B_{n}}
\cong R_{d}/(1-t_{1},\dots,1-t_{d}).\qedhere
\]
\end{proof}

Combining Lemmas~\ref{lem:TYM-tensor-decomposition} and~\ref{lem:coinvariants-TYMn-aDiag} with Proposition~\ref{prop:coinvariants-TYM}, we obtain the following result.

\begin{thm}\label{thm:coinvariants-TYM-tensor}
For $n\geq d+1$, there is an isomorphism of $R_{d}$-modules
\[
(V_{d})_{\B_{n}}
\cong
\bigoplus_{\{p_{1},\dots,p_{m}\}\in\cP_{d}}
\frac{R_{d}}{(1-t_{p_{1}},\dots,1-t_{p_{m}})}.
\]
Moreover, let $\bzeta=(\zeta_{1},\dots,\zeta_{d})\in(\bC^{\times})^{d}$, and write $\zeta_{p_{i}}:=\prod_{k\in p_{i}}\zeta_{k}$ for each $p_{i}\subseteq[d]$. Then
\[
(\TYM_{\zeta_{1}}\otimes\cdots\otimes
\TYM_{\zeta_{d}})_{\B_{n}}
\cong
\bigoplus_{\substack{
\{p_{1},\dots,p_{m}\}\in\cP_{d}\\
\zeta_{p_{1}}=\cdots=\zeta_{p_{m}}=1
}}
\bC.
\]
\end{thm}

\begin{proof}
By Lemma~\ref{lem:TYM-tensor-decomposition}, we have an isomorphism of $\B_{n}$-representations
\[
V_{d}
\cong
\bigoplus_{\{p_{1},\dots,p_{m}\}\in\cP_{d}}\aDiag\left(\TYM_{t_{p_{1}}}\otimes\cdots\otimes\TYM_{t_{p_{m}}}\right),
\]
where we have used the assumption $n\geq d+1$ to identify $\cP_{d}^{n}$ with $\cP_{d}$. Since taking coinvariants commutes with direct sums, we have
\[
(V_{d})_{\B_{n}}
\cong
\bigoplus_{\{p_{1},\dots,p_{m}\}\in\cP_{d}} \left(\aDiag\left(\TYM_{t_{p_{1}}}\otimes\cdots\otimes\TYM_{t_{p_{m}}}\right)\right)_{\B_{n}}.
\]
For partitions with $m\geq2$ parts, we apply Lemma~\ref{lem:coinvariants-TYMn-aDiag}, noting that $n\geq d+1\geq m+1$. For partitions with a single part, the corresponding summand is $\TYM_{t_{1}\cdots t_{d}}$, whose coinvariants are given by Proposition~\ref{prop:coinvariants-TYM}. In both cases, after extension of scalars to $R_{d}$, we obtain
\[
\left(
\aDiag\left(
\TYM_{t_{p_{1}}}\otimes\cdots\otimes
\TYM_{t_{p_{m}}}
\right)
\right)_{\B_{n}}
\cong
\frac{R_{d}}{(1-t_{p_{1}},\dots,1-t_{p_{m}})}.
\]
Combining these isomorphisms proves the first assertion.

We now consider the complex specialisations. Let $\bzeta=(\zeta_{1},\dots,\zeta_{d}) \in(\bC^{\times})^{d}$. Since taking coinvariants commutes with extension of scalars, we have
\[
(\TYM_{\zeta_{1}}\otimes\cdots\otimes
\TYM_{\zeta_{d}})_{\B_{n}}
\cong
\bigoplus_{\{p_{1},\dots,p_{m}\}\in\cP_{d}}
\frac{\bC}{(1-\zeta_{p_{1}},\dots,1-\zeta_{p_{m}})}.
\]
Since $\bC$ is a field, each summand is isomorphic to $\bC$ if $\zeta_{p_{1}}=\cdots=\zeta_{p_{m}}=1$, and vanishes otherwise. This gives the second assertion.
\end{proof}

Since $\TYM_{t}^{\vee}\cong\TYM_{t^{-1}}$ by~\eqref{eq:dual-TYM}, Theorem~\ref{thm:coinvariants-TYM-tensor} implies that $(V_{d}^{\vee})_{\B_{n}}$ is a torsion $R_{d}$-module. Together with Lemma~\ref{lem:invariants-coinvariants-duality}, this allows us to determine the invariants of $V_{d}$ and its complex specialisations.

\begin{coro}\label{coro:invariants-TYM-tensor}
For $n\geq d+1$, we have
\[
V_{d}^{\B_{n}}=0.
\]
Moreover, let $\bzeta=(\zeta_{1},\dots,\zeta_{d})\in(\bC^{\times})^{d}$, and write $\zeta_{p_{i}}:=\prod_{k\in p_{i}}\zeta_{k}$ for each $p_{i}\subseteq[d]$. Then
\[
(\TYM_{\zeta_{1}}\otimes\cdots\otimes\TYM_{\zeta_{d}})^{\B_{n}}
\cong
\bigoplus_{\substack{
\{p_{1},\dots,p_{m}\}\in\cP_{d}\\
\zeta_{p_{1}}=\cdots=\zeta_{p_{m}}=1
}}
\bC.
\]
\end{coro}

\begin{proof}
Since $\TYM_{t}^{\vee}\cong\TYM_{t^{-1}}$ by~\eqref{eq:dual-TYM}, we have $V_{d}^{\vee} \cong\TYM_{t_{1}^{-1}}\otimes\cdots\otimes \TYM_{t_{d}^{-1}}$. Applying Theorem~\ref{thm:coinvariants-TYM-tensor}, we obtain
\[
(V_{d}^{\vee})_{\B_{n}}
\cong
\bigoplus_{\{p_{1},\dots,p_{m}\}\in\cP_{d}}\frac{R_{d}}{(1-t_{p_{1}}^{-1},\dots,1-t_{p_{m}}^{-1})}.
\]
Since $R_{d}$ is an integral domain and each summand is a torsion $R_{d}$-module, so is $(V_{d}^{\vee})_{\B_{n}}$. Hence, by Lemma~\ref{lem:invariants-coinvariants-duality}, we have $V_{d}^{\B_{n}} \cong\mathrm{Hom}_{R_{d}}\bigl((V_{d}^{\vee})_{\B_{n}},R_{d}\bigr)=0$.

We now consider the complex specialisations. By~\eqref{eq:dual-TYM}, we have
\[
(\TYM_{\zeta_{1}}\otimes\cdots\otimes\TYM_{\zeta_{d}})^{\vee}
\cong
\TYM_{\zeta_{1}^{-1}}\otimes\cdots\otimes\TYM_{\zeta_{d}^{-1}}.
\]
Applying Theorem~\ref{thm:coinvariants-TYM-tensor} to these complex representations gives
\[
\bigl((\TYM_{\zeta_{1}}\otimes\cdots\otimes\TYM_{\zeta_{d}})^{\vee}\bigr)_{\B_{n}}
\cong
\bigoplus_{\substack{\{p_{1},\dots,p_{m}\}\in\cP_{d}\\\zeta_{p_{1}}^{-1}=\cdots=\zeta_{p_{m}}^{-1}=1}}\bC
\cong
\bigoplus_{\substack{\{p_{1},\dots,p_{m}\}\in\cP_{d}\\\zeta_{p_{1}}=\cdots=\zeta_{p_{m}}=1}} \bC.
\]
Finally, Lemma~\ref{lem:invariants-coinvariants-duality}, applied over $\bC$, identifies the invariants with the dual of these coinvariants. Since the latter form a finite-dimensional complex vector space, the result follows.
\end{proof}

\paragraph*{On the bound $n\geq d+1$.}
The assumption $n\geq d+1$ in Lemma~\ref{lem:coinvariants-TYMn-aDiag}, and then in Theorem~\ref{thm:coinvariants-TYM-tensor}, is essential to our computation. Indeed, the proof relies on the existence of an unoccupied index among $\{1,\dots,n\}$, which allows us to obtain the relations $(1-t_{k})v_{1,\dots,d}=0$ for every $1\leq k\leq d$. When $n=d$, no such index is available, and the resulting coinvariants may be different. For instance, when $n=d=2$, the anti-diagonal subrepresentation $\aDiag(\TYM_{t_{1}}\otimes\TYM_{t_{2}})$ is generated by $v_{1,2}$ and $v_{2,1}$. The defining relations of its coinvariants are $v_{1,2}=t_{1}v_{2,1}$ and $v_{2,1}=t_{2}v_{1,2}$, which yield
\[
\bigl(\aDiag(\TYM_{t_{1}}\otimes\TYM_{t_{2}})\bigr)_{\B_{2}}
\cong
R_{2}/(1-t_{1}t_{2}),
\]
rather than $R_{2}/(1-t_{1},1-t_{2})$, as obtained for $n\geq3$. Then, by Lemma~\ref{lem:TYM-tensor-decomposition}, we have
\[
(\TYM_{t_{1}}\otimes\TYM_{t_{2}})_{\B_{2}}
\cong
R_{2}/(1-t_{1}t_{2})
\oplus
R_{2}/(1-t_{1}t_{2}).
\]
Thus, the bound $n\geq d+1$ cannot be lowered in general. Although the decomposition of Lemma~\ref{lem:TYM-tensor-decomposition} remains valid for $n=d$, the computation of the corresponding anti-diagonal coinvariants requires a separate analysis, which we do not pursue here.

%--------------------------------------------------

\section{Other tensor products of low-dimensional representations}\label{s:other-tensor-products}

In this section, we compute the coinvariants of tensor products of two representations chosen among the Tong--Yang--Ma, unreduced Burau, and reduced Burau representations, together with their complex specialisations. The tensor product of two Tong--Yang--Ma representations has already been treated in \S\ref{ssec:tensor_tym}.

%--------------------------------------------------

\begin{prop}\label{prop:coinvariants-UBur2}
For $n\geq 2$, there is an isomorphism of $\bZ[t^{\pm1},q^{\pm1}]$-modules
\[
(\UBur_{t}\otimes\UBur_{q})_{\B_{n}}
\cong
\bZ[t^{\pm1},q^{\pm1}]
\oplus
\frac{\bZ[t^{\pm1},q^{\pm1}]}{(1-tq)}.
\]
Moreover, for $\bzeta=(\zeta_{t},\zeta_{q})\in(\bC^{\times})^{2}$, we have
\[
(\UBur_{\zeta_{t}}\otimes\UBur_{\zeta_{q}})_{\B_{n}}
\cong
\begin{cases}
\bC^{2} & \text{if }\zeta_{t}\zeta_{q}=1,\\
\bC & \text{otherwise}.
\end{cases}
\]
\end{prop}

\begin{proof}
Let $R:=\bZ[t^{\pm1},q^{\pm1}]$ and $\rho:=\UBur_{t}\otimes\UBur_{q}$. The underlying $R$-module $V_{\rho}$ has basis $\{v_{j,k}:=v_{j}\otimes v_{k} :1\leq j,k\leq n\}$. As in~\eqref{eq:def-invariants-coinvariants-Phi-Psi}, we regard the elements $(1-\rho(\sigma_{i}))(v_{j,k})$ for $1\leq i\leq n-1$ as the defining relations of $(V_{\rho})_{\B_{n}}$. First, we compute that $(1-\rho(\sigma_{j-1}))(v_{j,k})=v_{j,k}-v_{j-1,k}$ for $2\leq j<k\leq n$, while $(1-\rho(\sigma_{k-1}))(v_{1,k})=v_{1,k}-v_{1,k-1}$ for $3\leq k\leq n$. Therefore, $v_{j,k}=v_{1,2}$ in $(V_{\rho})_{\B_{n}}$ whenever $j<k$. Similarly, we have the relations $v_{j,k}-v_{j,k-1}=0$ for $2\leq k<j\leq n$, and $v_{j,1}-v_{j-1,1}=0$ for $3\leq j\leq n$. This implies that $v_{j,k}=v_{2,1}$ in $(V_{\rho})_{\B_{n}}$ whenever $j>k$. Finally, for $2\leq j\leq n$, we have $(1-\rho(\sigma_{j-1}))(v_{j,j})=v_{j,j}-v_{j-1,j-1}$, so all diagonal basis elements are identified with $v_{1,1}$. Thus, $(V_{\rho})_{\B_{n}}$ is generated by $v_{1,1}$, $v_{1,2}$, and $v_{2,1}$. Computing the remaining relations from the action of the standard braid generators, we obtain
\[
\begin{cases}
v_{1,2}-(1-t)v_{1,1}-tv_{2,1}=0,\\
v_{2,1}-(1-q)v_{1,1}-qv_{1,2}=0,\\
(t+q-2tq)v_{1,1}
-(1-t)qv_{1,2}
-t(1-q)v_{2,1}=0.
\end{cases}
\]
These relations are exhaustive: every other relation arising from the generators $\sigma_{i}$ becomes either trivial or one of the preceding relations under the identifications already established. The first relation gives $v_{1,2}=(1-t)v_{1,1}+tv_{2,1}$. Substituting this into the second relation yields $(1-tq)(v_{2,1}-v_{1,1})=0$. The third relation then reduces to $t(1-tq)(v_{1,1}-v_{2,1})=0$, which is redundant.
Setting $w:=v_{2,1}-v_{1,1}$, we conclude that $(V_{\rho})_{\B_{n}}$ is generated by $v_{1,1}$ and $w$, subject to the single relation $(1-tq)w=0$. Therefore, we have
\[
(\UBur_{t}\otimes\UBur_{q})_{\B_{n}}
\cong
\frac{R\{v_{1,1},w\}}{((1-tq)w)}
\cong
R\oplus R/(1-tq).
\]

We now consider the complex specialisations.
Let $\bzeta=(\zeta_{t},\zeta_{q})
\in(\bC^{\times})^{2}$.
Since taking coinvariants commutes with
extension of scalars, we have
\[
(\UBur_{\zeta_{t}}\otimes
\UBur_{\zeta_{q}})_{\B_{n}}
\cong
(\UBur_{t}\otimes\UBur_{q})_{\B_{n}}
\otimes_{R}\bC_{\bzeta}
\cong
\bC\oplus\bC/(1-\zeta_{t}\zeta_{q}).
\]
The second summand is isomorphic to $\bC$ if $\zeta_{t}\zeta_{q}=1$, and vanishes otherwise.
\end{proof}

\begin{prop}\label{prop:coinvariants-TYMxUBur}
For $n\geq 3$, there is an isomorphism of $\bZ[t^{\pm1},q^{\pm1}]$-modules
\[
(\TYM_{t}\otimes\UBur_{q})_{\B_{n}}
\cong
\frac{\bZ[t^{\pm1},q^{\pm1}]}{(1-t)} \oplus \frac{\bZ[t^{\pm1},q^{\pm1}]}{(1-t,1-q)}.
\]
Moreover, for $\bzeta=(\zeta_{t},\zeta_{q})\in(\bC^{\times})^{2}$, we have
\[
(\TYM_{\zeta_{t}}\otimes\UBur_{\zeta_{q}})_{\B_{n}}
\cong
\begin{cases}
\bC^{2}
& \text{if }\zeta_{t}=\zeta_{q}=1,\\
\bC
& \text{if }\zeta_{t}=1\text{ and }\zeta_{q}\neq 1,\\
0
& \text{if }\zeta_{t}\neq 1.
\end{cases}
\]
\end{prop}

\begin{proof}
Let $R:=\bZ[t^{\pm1},q^{\pm1}]$ and $\rho:=\TYM_{t}\otimes\UBur_{q}$. The underlying $R$-module $V_{\rho}$ has basis $\{v_{j,k}:=v_{j}\otimes v_{k}:1\leq j,k\leq n\}$. As in~\eqref{eq:def-invariants-coinvariants-Phi-Psi}, we regard the elements $(1-\rho(\sigma_{i}))(v_{j,k})$ as the defining relations of $(V_{\rho})_{\B_{n}}$. We first observe that all basis elements $v_{j,k}$ with $j<k$ are identified in the coinvariant module. First, we compute that $(1-\rho(\sigma_{j-1}))(v_{j,k}) =v_{j,k}-v_{j-1,k}$ for $2\leq j<k\leq n$, while we have $(1-\rho(\sigma_{k-1}))(v_{1,k}) =v_{1,k}-v_{1,k-1}$ for $3\leq k\leq n$. Hence, we have $v_{j,k}=v_{1,2}$ in $(V_{\rho})_{\B_{n}}$ whenever $j<k$. Similarly, we have $(1-\rho(\sigma_{k-1}))(v_{j,k}) =v_{j,k}-v_{j,k-1}$ for $2\leq k<j\leq n$, and $(1-\rho(\sigma_{j-1}))(v_{j,1}) =v_{j,1}-v_{j-1,1}$ for $3\leq j\leq n$. Hence, $v_{j,k}=v_{2,1}$ whenever $j>k$ in $(V_{\rho})_{\B_{n}}$. Furthermore, for $2\leq j\leq n$, we compute that $(1-\rho(\sigma_{j-1}))(v_{j,j}) =v_{j,j}-v_{j-1,j-1}$, so all diagonal basis elements are identified with $v_{1,1}$ in $(V_{\rho})_{\B_{n}}$. Finally, applying $1-\rho(\sigma_{1})$ to $v_{1,2}$ gives $v_{1,2}-tv_{2,1}=0$. Therefore, $(V_{\rho})_{\B_{n}}$ is generated by $v_{1,1}$ and $v_{2,1}$.

We now determine the relations between these generators. Since $n\geq3$, applying $1-\rho(\sigma_{2})$ to $v_{2,1}$ gives $(1-t)v_{2,1}=0$. Furthermore, applying $1-\rho(\sigma_{1})$ to $v_{2,1}$ and $v_{1,1}$, respectively, yields
\[
\begin{cases}
v_{2,1}-(1-q)v_{1,1}-qv_{1,2}=0,\\
v_{1,1}-t(1-q)v_{2,1}-tqv_{1,1}=0.
\end{cases}
\]
Using $v_{1,2}=tv_{2,1}$, we therefore obtain
\[
\begin{cases}
(1-t)v_{2,1}=0,\\
(1-tq)v_{2,1}-(1-q)v_{1,1}=0,\\
(1-tq)v_{1,1}-t(1-q)v_{2,1}=0.
\end{cases}
\]
We make the change of variable $w:=v_{2,1}-v_{1,1}$, so that $v_{1,1}=v_{2,1}-w$. The preceding relations become
\[
\begin{cases}
(1-t)v_{2,1}=0,\\
q(1-t)v_{2,1}+(1-q)w=0,\\
(1-t)v_{2,1}-(1-tq)w=0.
\end{cases}
\]
Using the first relation, the last two reduce to $(1-q)w=0$ and $(1-tq)w=0$. Since $1-t=(1-tq)-t(1-q)$, these are equivalent to $(1-t)w=(1-q)w=0$. Therefore, $(V_{\rho})_{\B_{n}}$ is generated by $v_{2,1}$ and $w$, satisfying $(1-t)v_{2,1}=0$ and $(1-t)w=(1-q)w=0$. In particular, there is a surjective $R$-module homomorphism
\[
\frac{R}{(1-t)}
\oplus
\frac{R}{(1-t,1-q)}
\longtwoheadrightarrow (V_{\rho})_{\B_{n}}
\]
sending the standard generators to $v_{2,1}$ and $w$, respectively. It remains to show that there are no additional relations. We define the $R$-linear maps $\varepsilon_{1}\colon V_{\rho} \to R/(1-t)$ and $\varepsilon_{2}\colon V_{\rho} \to R/(1-t,1-q)$ by $\varepsilon_{1}(v_{j,k})= 1$ and $\varepsilon_{2}(v_{j,k})= -\delta_{j,k}$ respectively. Modulo $(1-t)$, the first factor becomes the permutation representation of $\Sym_{n}$, and the augmentation of $\UBur_{q}$ is $\B_{n}$-invariant, so the map $\varepsilon_{1}$ is $\B_{n}$-invariant. Likewise, modulo $(1-t,1-q)$, both factors become permutation representations, so $\varepsilon_{2}$ is invariant since simultaneous permutations preserve the equality of indices. Therefore, the pair $(\varepsilon_{1},\varepsilon_{2})$ induces an $R$-module homomorphism
\[
(V_{\rho})_{\B_{n}}
\longrightarrow
\frac{R}{(1-t)}
\oplus
\frac{R}{(1-t,1-q)}.
\]
By construction, we have $(\varepsilon_{1},\varepsilon_{2})(v_{2,1})=(1,0)$ and $(\varepsilon_{1},\varepsilon_{2})(w)=(0,1)$. Hence, its composite with the preceding surjection is the identity. We conclude that
\[
(\TYM_{t}\otimes\UBur_{q})_{\B_{n}}
\cong
\frac{R}{(1-t)}
\oplus
\frac{R}{(1-t,1-q)}.
\]

Finally, since taking coinvariants commutes with extension of scalars, we obtain
\[
(\TYM_{\zeta_{t}}\otimes
\UBur_{\zeta_{q}})_{\B_{n}}
\cong
(\TYM_{t}\otimes\UBur_{q})_{\B_{n}}
\otimes_{R}\bC_{\bzeta}
\qquad\cong
\frac{\bC}{(1-\zeta_{t})}
\oplus
\frac{\bC}{(1-\zeta_{t},1-\zeta_{q})}.
\]
Since $\bC$ is a field, this module is two-dimensional if $\zeta_{t}=\zeta_{q}=1$, one-dimensional if $\zeta_{t}=1$ and $\zeta_{q}\neq1$, and zero if $\zeta_{t}\neq1$.
\end{proof}

\begin{prop}\label{prop:coinvariants-Bur2}
For $n\geq 3$, there is an isomorphism of $\bZ[t^{\pm1},q^{\pm1}]$-modules
\[
(\Bur_{t}\otimes\Bur_{q})_{\B_{n}}
\cong
\bZ[t^{\pm1},q^{\pm1}]/(1-tq).
\]
Moreover, for
$\bzeta=(\zeta_{t},\zeta_{q})\in(\bC^{\times})^{2}$, we have
\[
(\Bur_{\zeta_{t}}\otimes\Bur_{\zeta_{q}})_{\B_{n}}
\cong
\begin{cases}
\bC & \text{if }\zeta_{t}\zeta_{q}=1,\\
0 & \text{otherwise}.
\end{cases}
\]
\end{prop}

\begin{proof}
Let $R:=\bZ[t^{\pm1},q^{\pm1}]$ and $\rho:=\Bur_{t}\otimes\Bur_{q}$. The underlying $R$-module $V_{\rho}$ has basis $\{v_{j,k}:=v_{j}\otimes v_{k}:1\leq j,k\leq n-1\}$. As in~\eqref{eq:def-invariants-coinvariants-Phi-Psi}, we regard the elements $(1-\rho(\sigma_{i}))(v_{j,k})$ as the defining relations of $(V_{\rho})_{\B_{n}}$. First, we compute that $(1-\rho(\sigma_{k}))(v_{j,k-1})=-qv_{j,k}$ for $1\leq j<k-1\leq n-2$, while $(1-\rho(\sigma_{j}))(v_{j-1,k})=-tv_{j,k}$ for $1\leq k<j-1\leq n-2$. Since $t$ and $q$ are invertible, we deduce that $v_{j,k}=0$ in $(V_{\rho})_{\B_{n}}$ whenever $|j-k|\geq2$. Thus, only the diagonal elements and those immediately above and below the diagonal may remain nonzero. For each $1\leq j\leq n-2$, applying $1-\rho(\sigma_{j})$ to $v_{j+1,j+1}$, $v_{j,j+1}$, and $v_{j+1,j}$ gives
\[
\begin{cases}
v_{j,j}+v_{j,j+1}+v_{j+1,j}=0,\\
tv_{j,j}+(1+t)v_{j,j+1}=0,\\
qv_{j,j}+(1+q)v_{j+1,j}=0.
\end{cases}
\]
Subtracting $t$ times the first relation from the second yields $v_{j,j+1}=tv_{j+1,j}$. The first relation then gives $v_{j,j}=-(1+t)v_{j+1,j}$. Substituting this into the third relation, we obtain $(1-tq)v_{j+1,j}=0$. Furthermore, applying $1-\rho(\sigma_{j+1})$ to $v_{j,j+1}$ gives $tqv_{j+1,j+1}+(1+q)v_{j,j+1}=0$. Using $v_{j,j+1}=tv_{j+1,j}$, we deduce that
\[
v_{j+1,j+1}
=-(1+q^{-1})v_{j+1,j}
=-(1+t)v_{j+1,j},
\]
where the last equality follows from $(1-tq)v_{j+1,j}=0$. It remains to identify the classes $v_{j+1,j}$. For $1\leq j\leq n-3$, applying $1-\rho(\sigma_{j+1})$ to $v_{j+2,j}$ and $v_{j+2,j+1}$ gives
\[
\begin{cases}
v_{j+1,j}+qv_{j+1,j+1}+qv_{j+2,j+1}=0,\\
qv_{j+1,j+1}+(1+q)v_{j+2,j+1}=0.
\end{cases}
\]
Subtracting these relations yields $v_{j+1,j}=v_{j+2,j+1}$. Therefore, all the $v_{j+1,j}$ have the same class, which we denote by $v_{2,1}$. We have therefore shown that in $(V_{\rho})_{\B_{n}}$
\[
v_{j,k}=
\begin{cases}
-(1+t)v_{2,1} & j=k,\\
tv_{2,1} & k=j+1,\\
v_{2,1} & j=k+1,\\
0 & |j-k|\geq2,
\end{cases}
\]
and that $(1-tq)v_{2,1}=0$. In particular, there is a surjective $R$-module homomorphism
\[
R/(1-tq)\longtwoheadrightarrow (V_{\rho})_{\B_{n}}
\]
sending the class of $1$ to $v_{2,1}$. It remains to show that there are no additional relations. We consider the $R$-linear map $\varepsilon\colon V_{\rho}\to R/(1-tq)$ defined on the basis elements by
\[
\varepsilon(v_{j,k})=
\begin{cases}
-(1+t) & j=k,\\
t & k=j+1,\\
1 & j=k+1,\\
0 & |j-k|\geq2.
\end{cases}
\]
Let $M$ denote the $(n-1)\times(n-1)$ matrix with entries $M_{j,k}=\varepsilon(v_{j,k})$. A direct computation using the defining matrices of the reduced Burau representations shows that modulo $(1-tq)$
\[
\Bur_{t}(\sigma_{i})^{\tr}M\Bur_{q}(\sigma_{i})=M
\]
for all $1\leq i\leq n-1$. Since the $(j,k)$-entry of the left-hand side is precisely $\varepsilon(\rho(\sigma_{i})v_{j,k})$, this proves that $\varepsilon$ is $\B_{n}$-invariant. Thus, it induces an $R$-module homomorphism
\[
(V_{\rho})_{\B_{n}}\longrightarrow R/(1-tq).
\]
Since $\varepsilon(v_{2,1})=1$, its composite with the preceding surjection from $R/(1-tq)$ is the identity. We conclude that
\[
(\Bur_{t}\otimes\Bur_{q})_{\B_{n}}
\cong R/(1-tq).
\]

Finally, since taking coinvariants commutes with extension of scalars, we obtain
\[
\begin{aligned}
(\Bur_{\zeta_{t}}\otimes
\Bur_{\zeta_{q}})_{\B_{n}}
&\cong
(\Bur_{t}\otimes\Bur_{q})_{\B_{n}}
\otimes_{R}\bC_{\bzeta}\\
&\cong
\bC/(1-\zeta_{t}\zeta_{q}).
\end{aligned}
\]
This quotient is isomorphic to $\bC$ if $\zeta_{t}\zeta_{q}=1$, and vanishes otherwise.
\end{proof}

\begin{prop}\label{prop:coinvariants-TYMxBur}
For $n\geq 3$, there is an isomorphism of $\bZ[t^{\pm1},q^{\pm1}]$-modules
\[
(\TYM_{t}\otimes\Bur_{q})_{\B_{n}}
\cong
\bZ[t^{\pm1},q^{\pm1}]/(1-t,1-q).
\]
Moreover, for
$\bzeta=(\zeta_{t},\zeta_{q})\in(\bC^{\times})^{2}$, we have
\[
(\TYM_{\zeta_{t}}\otimes\Bur_{\zeta_{q}})_{\B_{n}}
\cong
\begin{cases}
\bC & \text{if }\zeta_{t}=\zeta_{q}=1,\\
0 & \text{otherwise}.
\end{cases}
\]
\end{prop}

\begin{proof}
Let $R:=\bZ[t^{\pm1},q^{\pm1}]$ and $\rho:=\TYM_{t}\otimes\Bur_{q}$. The underlying $R$-module $V_{\rho}$ has basis $\{v_{j,k}:=v_{j}\otimes v_{k}: 1\leq j\leq n,\ 1\leq k\leq n-1\}$. As in~\eqref{eq:def-invariants-coinvariants-Phi-Psi}, we regard the elements $(1-\rho(\sigma_{i}))(v_{j,k})$ as the defining relations of $(V_{\rho})_{\B_{n}}$. First, we compute that $(1-\rho(\sigma_{k}))(v_{j,k+1}) =-v_{j,k}$ for $1\leq k\leq n-2$ and $j\notin\{k,k+1\}$. Similarly, we have $(1-\rho(\sigma_{n-1}))(v_{j,n-2}) =-qv_{j,n-1}$ for $j\notin\{n-1,n\}$. Since $q$ is invertible, we deduce that $v_{j,k}=0$ in $(V_{\rho})_{\B_{n}}$ whenever $j\neq k,k+1$. Hence, only the classes $v_{k,k}$ and $v_{k+1,k}$ for $1\leq k\leq n-1$ may remain nonzero. For $1\leq k\leq n-2$, applying $1-\rho(\sigma_{k})$ to $v_{k+1,k+1}$, $v_{k,k+1}$, $v_{k+1,k}$, and $v_{k,k}$ gives
\[
\begin{cases}
v_{k+1,k+1}=v_{k,k},\\
v_{k+1,k+1}=-v_{k+1,k},\\
v_{k+1,k}=-qv_{k,k},\\
v_{k,k}=-tqv_{k+1,k}.
\end{cases}
\]
Therefore, all the $v_{k,k}$ have the same class $v_{1,1}$, and $(1-q)v_{1,1}=(1-t)v_{1,1}=0$. Moreover, $v_{k+1,k}=-v_{1,1}$ for $1\leq k\leq n-2$, while the relation $v_{n,n-1}=-qv_{n-1,n-1}$ gives $v_{n,n-1}=-v_{1,1}$ as well. Hence, $(V_{\rho})_{\B_{n}}$ is generated by $v_{1,1}$, subject to $(1-t)v_{1,1}=(1-q)v_{1,1}=0$. In particular, there is a surjective $R$-module homomorphism
\[
R/(1-t,1-q)\longtwoheadrightarrow (V_{\rho})_{\B_{n}}
\]
sending the class of $1$ to $v_{1,1}$. To show that there are no additional relations, we consider the $R$-linear map $\varepsilon\colon V_{\rho}\longrightarrow R/(1-t,1-q)$ defined by $v_{j,k}\mapsto \delta_{j,k}-\delta_{j,k+1}$. This map is $\B_{n}$-invariant, since at $t=q=1$ it coincides with the restriction of the standard $\Sym_{n}$-invariant pairing on the permutation representation to its tensor product with the standard subrepresentation, identified with the span of the vectors $e_{k}-e_{k+1}$. Hence, $\varepsilon$ factors through $(V_{\rho})_{\B_{n}}$, providing a left inverse to the preceding surjection, since $\varepsilon(v_{1,1})=1$. Therefore,
\[
(\TYM_{t}\otimes\Bur_{q})_{\B_{n}}
\cong R/(1-t,1-q).
\]

Finally, taking coinvariants commutes with extension of scalars, giving
\[
(\TYM_{\zeta_{t}}\otimes \Bur_{\zeta_{q}})_{\B_{n}}
\cong
(\TYM_{t}\otimes\Bur_{q})_{\B_{n}}
\otimes_{R}\bC_{\bzeta}\\
\cong
\bC/(1-\zeta_{t},1-\zeta_{q}).
\]
This is isomorphic to $\bC$ if $\zeta_{t}=\zeta_{q}=1$, and vanishes otherwise.
\end{proof}

\begin{prop}\label{prop:coinvariants-UBurxBur}
For $n\geq 3$, there is an isomorphism of $\bZ[t^{\pm1},q^{\pm1}]$-modules
\[
(\UBur_{t}\otimes\Bur_{q})_{\B_{n}}
\cong
\bZ[t^{\pm1},q^{\pm1}]/(1-tq).
\]
Moreover, for $\bzeta=(\zeta_{t},\zeta_{q})\in(\bC^{\times})^{2}$, we have
\[
(\UBur_{\zeta_{t}}\otimes\Bur_{\zeta_{q}})_{\B_{n}}
\cong
\begin{cases}
\bC & \text{if }\zeta_{t}\zeta_{q}=1,\\
0 & \text{otherwise}.
\end{cases}
\]
\end{prop}

\begin{proof}
Let $R:=\bZ[t^{\pm1},q^{\pm1}]$ and $\rho:=\UBur_{t}\otimes\Bur_{q}$. The underlying $R$-module $V_{\rho}$ has basis $\{v_{j,k}:=v_{j}\otimes v_{k}:1\leq j\leq n,\ 1\leq k\leq n-1\}$. As in~\eqref{eq:def-invariants-coinvariants-Phi-Psi}, we regard the elements $(1-\rho(\sigma_{i}))(v_{j,k})$ as the defining relations of $(V_{\rho})_{\B_{n}}$. First, we compute that $(1-\rho(\sigma_{k}))(v_{j,k-1})=-qv_{j,k}$ for $1\leq j<k\leq n-1$. Since $q$ is invertible, we deduce that $v_{j,k}=0$ in $(V_{\rho})_{\B_{n}}$ whenever $j<k$. Similarly, we have $(1-\rho(\sigma_{k}))(v_{k+2,k+1})=-v_{k+2,k}$ for $1\leq k\leq n-2$, so $v_{k+2,k}=0$ in $(V_{\rho})_{\B_{n}}$. Moreover, we have $(1-\rho(\sigma_{j-1}))(v_{j,k})=v_{j,k}-v_{j-1,k}$ for $j\geq k+3$. Hence $v_{j,k}=0$ whenever $j\geq k+2$ in $(V_{\rho})_{\B_{n}}$. Therefore, only the diagonal elements $v_{j,j}$ and the elements $v_{j+1,j}$ immediately below the diagonal may remain nonzero. For $2\leq j\leq n-1$, we have
\[
(1-\rho(\sigma_{j-1}))(v_{j,j}) =v_{j,j}-v_{j-1,j-1}-v_{j-1,j}.
\]
Since $v_{j-1,j}=0$, all diagonal elements have the same class in $(V_{\rho})_{\B_{n}}$. Furthermore, for $1\leq j\leq n-1$, we have
\[
(1-\rho(\sigma_{j}))(v_{j,j}) =(1+q-tq)v_{j,j}+tqv_{j+1,j}.
\]
Since $tq$ is invertible, all the elements $v_{j+1,j}$ also have the same class $v_{2,1}$ in $(V_{\rho})_{\B_{n}}$. Hence $(V_{\rho})_{\B_{n}}$ is generated by $v_{1,1}$ and $v_{2,1}$.

We now determine the relations between these generators. Applying $1-\rho(\sigma_{1})$ to $v_{2,1}$ gives $v_{2,1}+qv_{1,1}=0$, so $v_{2,1}=-qv_{1,1}$. On the other hand, applying $1-\rho(\sigma_{1})$ to $v_{1,2}$ gives $tv_{1,2}-(1-t)v_{1,1} -tv_{2,1}-tv_{2,2}=0$. Using $v_{1,2}=0$, $v_{2,2}=v_{1,1}$, and $v_{2,1}=-qv_{1,1}$, this reduces to $(1-tq)v_{1,1}=0$. Therefore, $(V_{\rho})_{\B_{n}}$ is generated by $v_{1,1}$, subject to $(1-tq)v_{1,1}=0$. In particular, there is a surjective $R$-module homomorphism
\[
R/(1-tq)\longtwoheadrightarrow(V_{\rho})_{\B_{n}}
\]
sending the class of $1$ to $v_{1,1}$. It remains to show that there are no additional relations. We consider the $R$-linear map $\varepsilon\colon V_{\rho}\to R/(1-tq)$ defined on the basis elements by $\varepsilon(v_{j,k})=\delta_{j,k}-q\delta_{j,k+1}$. Let $M$ be the $n\times(n-1)$ matrix with entries $M_{j,k}=\delta_{j,k}-q\delta_{j,k+1}$. A direct computation using the defining matrices of the unreduced and reduced Burau representations gives
\[
\UBur_{t}(\sigma_{i})^{\tr} M\Bur_{q}(\sigma_{i})-M
=(tq-1)e_{i}\bigl(-qf_{i-1}+(1+q)f_{i}-f_{i+1}\bigr)^{\tr},
\]
where $e_{i}$ and $f_{i}$ denote the standard basis vectors of $R^{n}$ and $R^{n-1}$, respectively, with $f_{0}=f_{n}=0$. Since $tq=1$ in $R/(1-tq)$, this proves that $\varepsilon$ is $\B_{n}$-invariant. Therefore, $\varepsilon$ induces an $R$-module homomorphism
\[
(V_{\rho})_{\B_{n}}\longrightarrow R/(1-tq).
\]
Since $\varepsilon(v_{1,1})=1$, its composite with the preceding surjection from $S$ is the identity. We conclude that
\[
(\UBur_{t}\otimes\Bur_{q})_{\B_{n}}
\cong R/(1-tq).
\]

We now consider the complex specialisations. Let $\bzeta=(\zeta_{t},\zeta_{q})\in(\bC^{\times})^{2}$. Since taking coinvariants commutes with extension of scalars, we obtain
\[
(\UBur_{\zeta_{t}}\otimes\Bur_{\zeta_{q}})_{\B_{n}}
\cong
(\UBur_{t}\otimes\Bur_{q})_{\B_{n}}
\otimes_{R}\bC_{\bzeta}
\cong
\bC/(1-\zeta_{t}\zeta_{q}).
\]
Since $\bC$ is a field, this quotient is one-dimensional if $\zeta_{t}\zeta_{q}=1$, and vanishes otherwise.
\end{proof}

\phantomsection
\addcontentsline{toc}{section}{References}
\renewcommand{\bibfont}{\normalfont\small}
\setlength{\bibitemsep}{0pt}
\printbibliography

@article {RWW,
	AUTHOR = {Randal-Williams, Oscar and Wahl, Nathalie},
	TITLE = {Homological stability for automorphism groups},
	JOURNAL = {Adv. Math.},
	FJOURNAL = {Advances in Mathematics},
	VOLUME = {318},
	YEAR = {2017},
	PAGES = {534--626},
	ISSN = {0001-8708},
	MRCLASS = {20J05},
	MRNUMBER = {3689750},
	URL = {https://doi.org/10.1016/j.aim.2017.07.022},
	Addendum = {\{\href{http://arxiv.org/abs/1409.3541}{arxiv:1409.3541}\}},
}

@article {S1,
    AUTHOR = {Souli\'{e}, Arthur},
     TITLE = {The {L}ong-{M}oody construction and polynomial functors},
   JOURNAL = {Ann. Inst. Fourier (Grenoble)},
  FJOURNAL = {Universit\'{e} de Grenoble. Annales de l'Institut Fourier},
    VOLUME = {69},
      YEAR = {2019},
    NUMBER = {4},
     PAGES = {1799--1856},
      ISSN = {0373-0956},
   MRCLASS = {18A25 (18M15 20C99 20F36 20J99)},
  MRNUMBER = {4010870},
MRREVIEWER = {Jingyan Li},
       URL = {http://aif.cedram.org/item?id=AIF_2019__69_4_1799_0},
}

@article {KrammerLK,
    AUTHOR = {Krammer, Daan},
     TITLE = {Braid groups are linear},
   JOURNAL = {Ann. of Math. (2)},
  FJOURNAL = {Annals of Mathematics. Second Series},
    VOLUME = {155},
      YEAR = {2002},
    NUMBER = {1},
     PAGES = {131--156},
      ISSN = {0003-486X},
   MRCLASS = {20F36 (57M07)},
  MRNUMBER = {1888796},
MRREVIEWER = {Luis Paris},
       DOI = {10.2307/3062152},
       URL = {http://dx.doi.org/10.2307/3062152},
}

@article {bigelow2001braid,
    AUTHOR = {Bigelow, Stephen J.},
     TITLE = {Braid groups are linear},
   JOURNAL = {J. Amer. Math. Soc.},
  FJOURNAL = {Journal of the American Mathematical Society},
    VOLUME = {14},
      YEAR = {2001},
    NUMBER = {2},
     PAGES = {471--486},
      ISSN = {0894-0347},
   MRCLASS = {20F36 (57M07)},
  MRNUMBER = {1815219},
MRREVIEWER = {Darryl McCullough},
       DOI = {10.1090/S0894-0347-00-00361-1},
       URL = {https://doi.org/10.1090/S0894-0347-00-00361-1},
}

@incollection {bigelowthelkrepresentation,
    AUTHOR = {Bigelow, Stephen},
     TITLE = {The {L}awrence-{K}rammer representation},
 BOOKTITLE = {Topology and geometry of manifolds ({A}thens, {GA}, 2001)},
    SERIES = {Proc. Sympos. Pure Math.},
    VOLUME = {71},
     PAGES = {51--68},
 PUBLISHER = {Amer. Math. Soc., Providence, RI},
      YEAR = {2003},
   MRCLASS = {20F36 (20C08)},
  MRNUMBER = {2024629},
MRREVIEWER = {Luisa Paoluzzi},
       DOI = {10.1090/pspum/071/2024629},
       URL = {http://dx.doi.org/10.1090/pspum/071/2024629},
}

@article {Lawrence1,
    AUTHOR = {Lawrence, R. J.},
     TITLE = {Homological representations of the {H}ecke algebra},
   JOURNAL = {Comm. Math. Phys.},
  FJOURNAL = {Communications in Mathematical Physics},
    VOLUME = {135},
      YEAR = {1990},
    NUMBER = {1},
     PAGES = {141--191},
      ISSN = {0010-3616},
   MRCLASS = {16G99 (20C32 20F36 32S40 57M25)},
  MRNUMBER = {1086755},
MRREVIEWER = {D. B. Fuchs},
       URL = {http://projecteuclid.org/euclid.cmp/1104201923},
}

@article{BirmanBrendlesurvey,
  title={Braids: a survey},
  author={Birman, Joan S. and Brendle, Tara E.},
  journal={Handbook of knot theory},
  pages={19--103},
  year={2005}
}

@inproceedings{burau,
  title={{\"U}ber Zopfgruppen und gleichsinnig verdrillte Verkettungen},
  author={Burau, Werner},
  booktitle={Abhandlungen aus dem Mathematischen Seminar der Universit{\"a}t Hamburg},
  volume={11},
  number={1},
  pages={179--186},
  year={1935},
  organization={Springer}
}

@article {Formanek,
    AUTHOR = {Formanek, Edward},
     TITLE = {Braid group representations of low degree},
   JOURNAL = {Proc. London Math. Soc. (3)},
  FJOURNAL = {Proceedings of the London Mathematical Society. Third Series},
    VOLUME = {73},
      YEAR = {1996},
    NUMBER = {2},
     PAGES = {279--322},
      ISSN = {0024-6115},
   MRCLASS = {20F36 (20C15 57M07)},
  MRNUMBER = {1397691},
MRREVIEWER = {J. S. Birman},
       DOI = {10.1112/plms/s3-73.2.279},
       URL = {https://doi-org.ezproxy.lib.gla.ac.uk/10.1112/plms/s3-73.2.279},
}

@article {Sysoeva1,
    AUTHOR = {Sysoeva, Inna},
     TITLE = {Dimension {$n$} representations of the braid group on {$n$}
              strings},
   JOURNAL = {J. Algebra},
  FJOURNAL = {Journal of Algebra},
    VOLUME = {243},
      YEAR = {2001},
    NUMBER = {2},
     PAGES = {518--538},
      ISSN = {0021-8693},
   MRCLASS = {20F36 (20C07)},
  MRNUMBER = {1850645},
MRREVIEWER = {Vitaly A. Roman\cprime kov},
       DOI = {10.1006/jabr.2001.8879},
       URL = {https://doi-org.ezproxy.lib.gla.ac.uk/10.1006/jabr.2001.8879},
}

@article {Sysoeva2,
    AUTHOR = {Sysoeva, Inna},
     TITLE = {Irreducible representations of braid group {$B_n$} of
              dimension {$n+1$}},
   JOURNAL = {J. Group Theory},
  FJOURNAL = {Journal of Group Theory},
    VOLUME = {24},
      YEAR = {2021},
    NUMBER = {1},
     PAGES = {39--78},
      ISSN = {1433-5883},
   MRCLASS = {20F36},
  MRNUMBER = {4193497},
MRREVIEWER = {Shoumin Liu},
       DOI = {10.1515/jgth-2020-0103},
       URL = {https://doi.org/10.1515/jgth-2020-0103},
}

@article {Formanek-Lee-Sysoeva-Vazirani,
    AUTHOR = {Formanek, Edward and Lee, Woo and Sysoeva, Inna and Vazirani,
              Monica},
     TITLE = {The irreducible complex representations of the braid group on
              {$n$} strings of degree {$\le n$}},
   JOURNAL = {J. Algebra Appl.},
  FJOURNAL = {Journal of Algebra and its Applications},
    VOLUME = {2},
      YEAR = {2003},
    NUMBER = {3},
     PAGES = {317--333},
      ISSN = {0219-4988},
   MRCLASS = {20F36 (20C08)},
  MRNUMBER = {1997750},
MRREVIEWER = {Graham C. Denham},
       DOI = {10.1142/S0219498803000477},
       URL = {https://doi-org.ezproxy.lib.gla.ac.uk/10.1142/S0219498803000477},
}

@article {ErdmannNakano,
    AUTHOR = {Erdmann, Karin and Nakano, Daniel K.},
     TITLE = {Representation type of {H}ecke algebras of type {$A$}},
   JOURNAL = {Trans. Amer. Math. Soc.},
  FJOURNAL = {Transactions of the American Mathematical Society},
    VOLUME = {354},
      YEAR = {2002},
    NUMBER = {1},
     PAGES = {275--285},
      ISSN = {0002-9947},
   MRCLASS = {20C08 (16G60)},
  MRNUMBER = {1859276},
MRREVIEWER = {Hiroshi Naruse},
       DOI = {10.1090/S0002-9947-01-02848-3},
       URL = {https://doi-org.ezproxy.lib.gla.ac.uk/10.1090/S0002-9947-01-02848-3},
}

@book {Brown,
    AUTHOR = {Brown, Kenneth S.},
     TITLE = {Cohomology of groups},
    SERIES = {Graduate Texts in Mathematics},
    VOLUME = {87},
      NOTE = {Corrected reprint of the 1982 original},
 PUBLISHER = {Springer-Verlag, New York},
      YEAR = {1994},
     PAGES = {x+306},
      ISBN = {0-387-90688-6},
   MRCLASS = {20J05 (20-02)},
  MRNUMBER = {1324339},
}

@article {Chen,
    AUTHOR = {Chen, Weiyan},
     TITLE = {Homology of braid groups, the {B}urau representation, and
              points on superelliptic curves over finite fields},
   JOURNAL = {Israel J. Math.},
  FJOURNAL = {Israel Journal of Mathematics},
    VOLUME = {220},
      YEAR = {2017},
    NUMBER = {2},
     PAGES = {739--762},
      ISSN = {0021-2172,1565-8511},
   MRCLASS = {20F36 (11G05 14H20 14H50)},
  MRNUMBER = {3666443},
MRREVIEWER = {Markus\ Szymik},
       DOI = {10.1007/s11856-017-1534-7},
       URL = {https://doi.org/10.1007/s11856-017-1534-7},
}

@article {CallegaroMoroniSalvetti,
    AUTHOR = {Callegaro, Filippo and Moroni, Davide and Salvetti, Mario},
     TITLE = {Cohomology of affine {A}rtin groups and applications},
   JOURNAL = {Trans. Amer. Math. Soc.},
  FJOURNAL = {Transactions of the American Mathematical Society},
    VOLUME = {360},
      YEAR = {2008},
    NUMBER = {8},
     PAGES = {4169--4188},
      ISSN = {0002-9947,1088-6850},
   MRCLASS = {20J06 (20F36)},
  MRNUMBER = {2395168},
MRREVIEWER = {Peter\ A.\ Linnell},
       DOI = {10.1090/S0002-9947-08-04488-7},
       URL = {https://doi.org/10.1090/S0002-9947-08-04488-7},
}

@book {KasselTuraev,
    AUTHOR = {Kassel, Christian and Turaev, Vladimir},
     TITLE = {Braid groups},
    SERIES = {Graduate Texts in Mathematics},
    VOLUME = {247},
      NOTE = {With the graphical assistance of Olivier Dodane},
 PUBLISHER = {Springer, New York},
      YEAR = {2008},
     PAGES = {xii+340},
      ISBN = {978-0-387-33841-5},
   MRCLASS = {20F36 (20C08 20F10 20F60 20M05 55R80 57M25 57M50)},
  MRNUMBER = {2435235},
MRREVIEWER = {Stephen\ P.\ Humphries},
       DOI = {10.1007/978-0-387-68548-9},
       URL = {https://doi.org/10.1007/978-0-387-68548-9},
}

@article {TYM,
    AUTHOR = {Tong, Dian-Min and Yang, Shan-De and Ma, Zhong-Qi},
     TITLE = {A new class of representations of braid groups},
   JOURNAL = {Comm. Theoret. Phys.},
  FJOURNAL = {Communications in Theoretical Physics},
    VOLUME = {26},
      YEAR = {1996},
    NUMBER = {4},
     PAGES = {483--486},
      ISSN = {0253-6102},
   MRCLASS = {20F36 (20C15)},
  MRNUMBER = {1456851},
       DOI = {10.1088/0253-6102/26/4/483},
       URL = {https://doi.org/10.1088/0253-6102/26/4/483},
}

@incollection {BelitskiiSergeichuk,
    AUTHOR = {Belitskii, Genrich R. and Sergeichuk, Vladimir V.},
     TITLE = {Complexity of matrix problems},
      NOTE = {Ninth Conference of the International Linear Algebra Society
              (Haifa, 2001)},
   JOURNAL = {Linear Algebra Appl.},
  FJOURNAL = {Linear Algebra and its Applications},
    VOLUME = {361},
      YEAR = {2003},
     PAGES = {203--222},
      ISSN = {0024-3795},
   MRCLASS = {16G20 (15A21 16G60)},
  MRNUMBER = {1955562},
MRREVIEWER = {Diane Castonguay},
       DOI = {10.1016/S0024-3795(02)00391-9},
       URL = {https://doi.org/10.1016/S0024-3795(02)00391-9},
}

@preamble{
   "\def\cprime{$'$} "
}

@article {Levaillant-invariants,
    AUTHOR = {Levaillant, Claire},
     TITLE = {Invariant subspaces of the {L}awrence-{K}rammer
              representation},
   JOURNAL = {Proc. Amer. Math. Soc.},
  FJOURNAL = {Proceedings of the American Mathematical Society},
    VOLUME = {140},
      YEAR = {2012},
    NUMBER = {8},
     PAGES = {2599--2612},
      ISSN = {0002-9939,1088-6826},
   MRCLASS = {20F36 (20C08)},
  MRNUMBER = {2910748},
MRREVIEWER = {Craig\ J.\ Dodge},
       DOI = {10.1090/S0002-9939-2011-11107-9},
       URL = {https://doi.org/10.1090/S0002-9939-2011-11107-9},
}

@article {Levaillant-Irred,
    AUTHOR = {Levaillant, Claire},
     TITLE = {Irreducibility of the {L}awrence-{K}rammer representation of
              the {BMW} algebra of type {$A_{n-1}$}},
   JOURNAL = {C. R. Math. Acad. Sci. Paris},
  FJOURNAL = {Comptes Rendus Math\'{e}matique. Acad\'{e}mie des Sciences.
              Paris},
    VOLUME = {347},
      YEAR = {2009},
    NUMBER = {1-2},
     PAGES = {15--20},
      ISSN = {1631-073X,1778-3569},
   MRCLASS = {20F36},
  MRNUMBER = {2536742},
MRREVIEWER = {Stephen\ P.\ Humphries},
       DOI = {10.1016/j.crma.2008.11.011},
       URL = {https://doi.org/10.1016/j.crma.2008.11.011},
}

@article {LevaillantWales,
    AUTHOR = {Levaillant, Claire and Wales, David},
     TITLE = {Parameters for which the {L}awrence-{K}rammer representation
              is reducible},
   JOURNAL = {J. Algebra},
  FJOURNAL = {Journal of Algebra},
    VOLUME = {323},
      YEAR = {2010},
    NUMBER = {7},
     PAGES = {1966--1982},
      ISSN = {0021-8693,1090-266X},
   MRCLASS = {20F36},
  MRNUMBER = {2594657},
MRREVIEWER = {Eric\ C.\ Rowell},
       DOI = {10.1016/j.jalgebra.2009.12.021},
       URL = {https://doi.org/10.1016/j.jalgebra.2009.12.021},
}

@article {Zinno,
    AUTHOR = {Zinno, Matthew G.},
     TITLE = {On {K}rammer's representation of the braid group},
   JOURNAL = {Math. Ann.},
  FJOURNAL = {Mathematische Annalen},
    VOLUME = {321},
      YEAR = {2001},
    NUMBER = {1},
     PAGES = {197--211},
      ISSN = {0025-5831,1432-1807},
   MRCLASS = {20F36 (57M25)},
  MRNUMBER = {1857374},
MRREVIEWER = {Vitaly\ A.\ Roman\cprime kov},
       DOI = {10.1007/PL00004501},
       URL = {https://doi.org/10.1007/PL00004501},
}

@misc{AnnexLS1,
    title        = {Annex file to the article: “Twisted homological stability for handlebody mapping class groups”},
    Author = {Lindell, Erik and Souli{\'e}, Arthur},
    howpublished = {Available \href{https://drive.google.com/file/d/13E8QrVAup85L_wgNxIWXYLrWKSAiMFUr/view}{here}},
    year         = {2025},
}

\noindent {Anh Trong Nam Hoang, \itshape Northeastern University, Boston, MA 02115, USA.}\\
\noindent {Email address: \tt an.hoang@northeastern.edu}

\noindent {Arthur Souli{\'e}, \itshape CNRS, LMNO, Normandie University, UNICAEN, 14000 Caen, France.}\\
\noindent {Email addresses: \tt artsou@hotmail.fr, arthur.soulie@unicaen.fr, arthur.soulie@cnrs.fr.}

\end{document}